\documentclass[preprint]{elsarticle}       
\usepackage{lipsum}
\makeatletter
\def\ps@pprintTitle{%
 \let\@oddhead\@empty
 \let\@evenhead\@empty
 \def\@oddfoot{}%
 \let\@evenfoot\@oddfoot}
\makeatother

\usepackage{lineno,hyperref}
\usepackage{rotating}
\usepackage{makeidx}
\usepackage{float}
\usepackage{epsfig}
\usepackage{amsmath,amssymb}
\usepackage[latin1]{inputenc}
\usepackage{epic,eepic}
\usepackage{overpic}
\usepackage{color}
\usepackage{amsfonts}
\usepackage{graphicx}
\usepackage{mwe} 
\usepackage{subcaption}

\usepackage{algorithm}
\usepackage{algorithmicx}
\usepackage{algpseudocode}
\usepackage{enumitem}

            {\hfill\penalty10000\raisebox{-.09em}{\large\bf\rm $\blacksquare$}\par\medskip}

\newtheorem{theorem}{Theorem}[section]
\newtheorem{lemma}[theorem]{Lemma}
\newtheorem{proposition}[theorem]{Proposition}

\newtheorem{definition}{Definition}

\DeclareMathOperator*{\argmin}{arg\,min}

\newcommand{\xx}{\mathbf{x}}

\journal{}

\begin{document}

\small

\begin{frontmatter}

\title{Non-linear Partition of Unity method}\tnotetext[label1]{Dionisio F. Y\'a\~nez has been supported through project CIAICO/2021/227 (Proyecto financiado por la Conselleria de Innovaci\'on, Universidades, Ciencia y Sociedad digital de la Generalitat Valenciana), by grant PID2020-117211GB-I00 and by PID2023-146836NB-I00 funded by MCIN/AEI/10.13039/501100011033.}

\author[UV]{Jos\'e M. Ram\'on}
\ead{Jose.Manuel.Ramon@uv.es}
\author[UPCT]{Juan Ruiz-\'Alvarez}
\ead{juan.ruiz@upct.es}
\author[UV]{Dionisio F. Y\'a\~nez}
\ead{Dionisio.Yanez@uv.es}

\date{Received: date / Accepted: date}

\address[UV]{Departamento de Matem\'aticas. Universidad de Valencia, Valencia (Spain).}
\address[UPCT]{Departamento de Matem\'atica Aplicada y Estad\'istica. Universidad  Polit\'ecnica de Cartagena, Cartagena (Spain).}


\begin{abstract}
This paper introduces the Non-linear Partition of Unity Method, a novel technique integrating Radial Basis Function interpolation and Weighted Essentially Non-Oscillatory algorithms. It addresses challenges in high-accuracy approximations, particularly near discontinuities, by adapting weights dynamically. The method is rooted in the Partition of Unity framework, enabling efficient decomposition of large datasets into subproblems while maintaining accuracy. Smoothness indicators and compactly supported functions ensure precision in regions with discontinuities. Error bounds are calculated and validate its effectiveness, showing improved interpolation in discontinuous and smooth regions. Some numerical experiments are performed to check the theoretical results.
\end{abstract}
\end{frontmatter}

\section{Introduction}\label{introduction}
Kernel-based methods are useful and powerful tools employed in several fields of science as computer aided design, numerical resolution of partial differential equations, image processing or others. These schemes are satisfactory used in interpolation, regression and machine learning due to their easy and quick implementation and their simple generalization for any dimension.

\begin{figure}
\begin{center}
		\begin{tabular}{cc}
			 \hspace{-1cm}	\includegraphics[width=8cm]{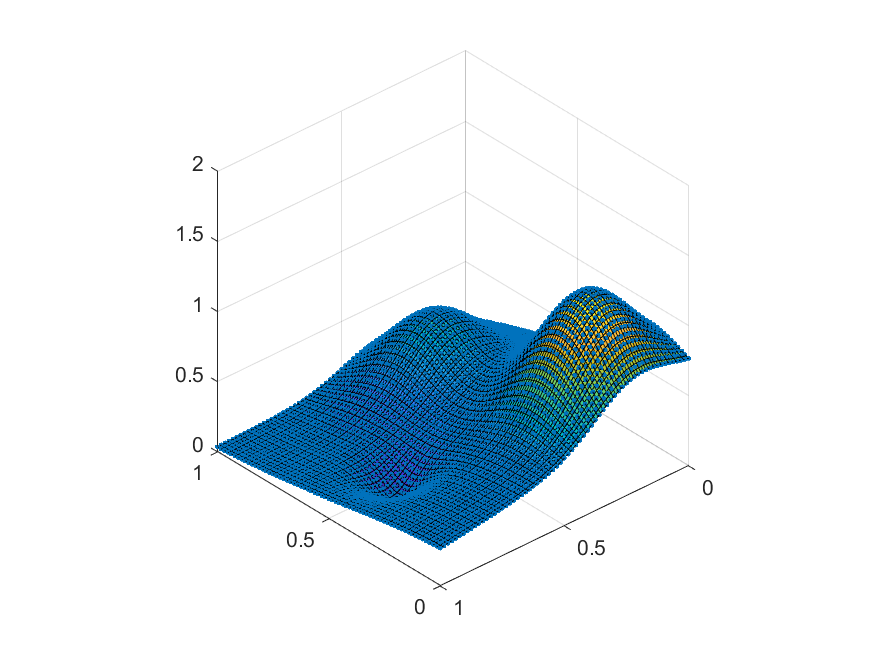} &
		 \hspace{-1cm}	\includegraphics[width=8cm]{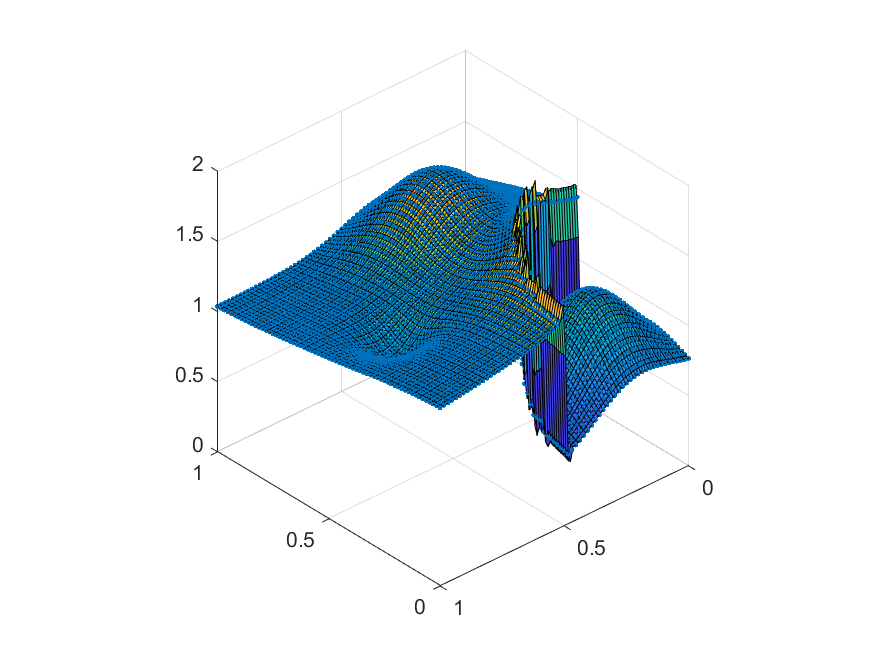}\\ (a) & (b)\\
  \hspace{-1cm}\includegraphics[width=6cm]{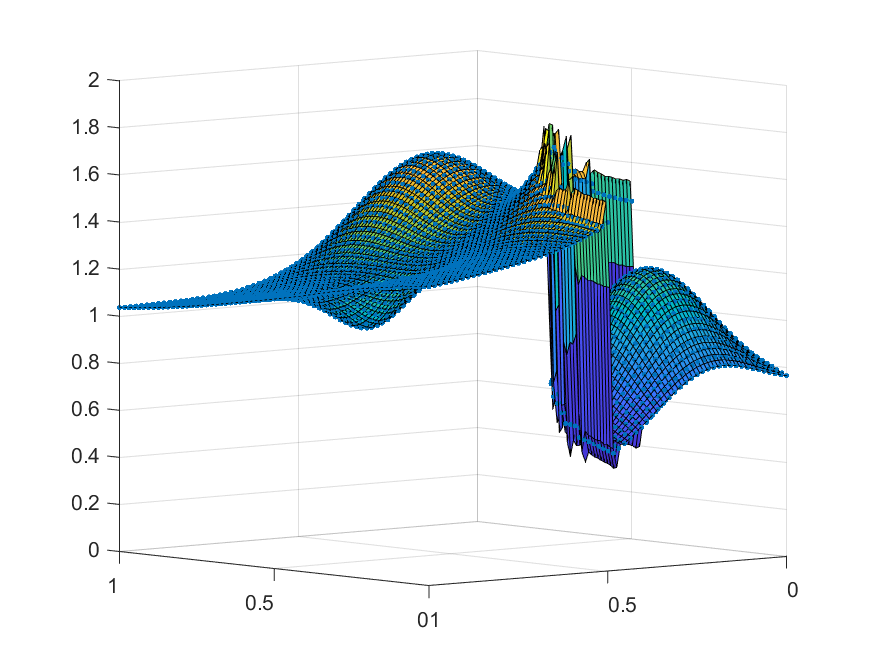} &
  \hspace{-1cm} \includegraphics[width=8cm]{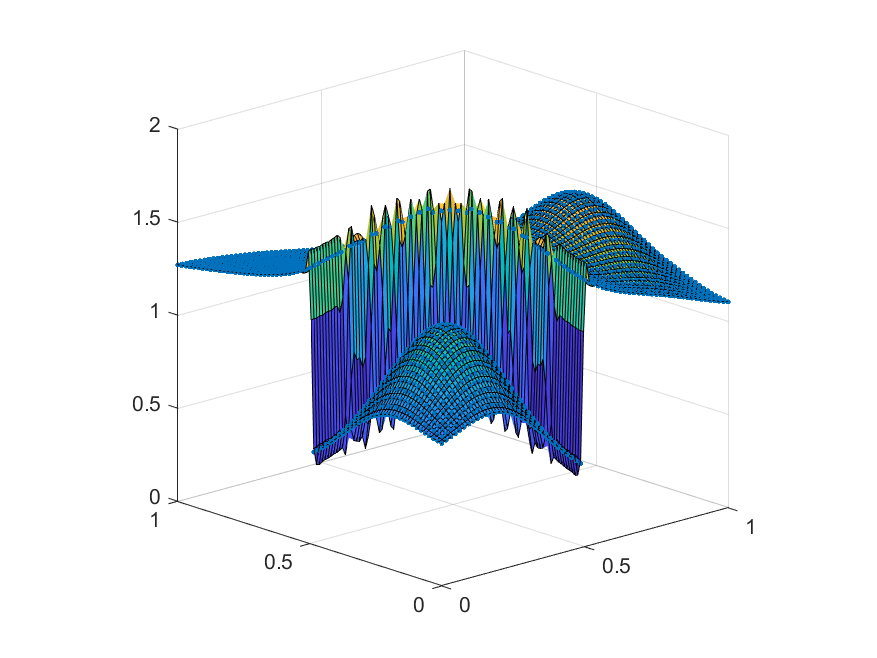}\\(c) & (d)
		\end{tabular}
\end{center}
\caption{(a) Approximation to Franke's function using PUM, (b), (c), (d) Approximation to function $f_1$, Eq. \eqref{frankesdisc} using PUM.}
	\label{figureintro}
	\end{figure}

The problem consists in finding a good approximation to an unknown function, $f:\mathbb{R}^n\to \mathbb{R}$, from some data values $F=\{f_i=f(\xx_i),\,\,i=1,\hdots,N\}$ given a set of arbitrarily distributed points on an open and bounded domain $\Omega \subset \mathbb{R}^n$, $X=\{\xx_i\in \Omega,\,\,i=1,\hdots,N\}$. Radial basis function method (RBF) consists on developing the function as a linear combination of a basis of radial functions centered in the data points $X$. There exists a vast literature about it (see, e.g. \cite{buhmann,FASSHAUER,wendland}), however the efficiency of the method can be improved in different ways. One approximation, called {\it partition of unity method} (PUM or PU method) (see e.g. \cite{buhmann,cavo15,cavo16,FASSHAUER,wendland}) involves partitioning the entire domain into multiple smaller subdomains and applying the RBF method to each of them. Ultimately, a convex combination of these interpolants yields the final approximation. Satisfactory results are obtained when these methods are employed to data originating from continuous function, see for example, the approximation presented in Fig. \ref{figureintro} (a) where we approximate the Franke's function:
\begin{equation}\label{frankesfunction}
\begin{split}
f(x,y)=&\frac34 e^{-1/4((9x-2)^2+(9y-2)^2)}+\frac34 e^{-1/49(9x+1)^2-1/10(9y+1)}\\
&+\frac12 e^{-1/4((9x-7)^2+(9y-3)^2)}-\frac15 e^{-(9x-4)^2-(9y-7)^2}
\end{split}
\end{equation}
using the PU method with data points located at $X=\{(i/64,j/64):i,j=0,\hdots,64\}$. However, these methods yield inadequate results when the data originates from functions with discontinuities. See for example Figs. \ref{figureintro} (b), (c) and (d), where we present the result obtained when applying the PUM to the piecewise smooth function:
\begin{equation}\label{frankesdisc}
f_1(x,y)=\left\{
                \begin{array}{ll}
                  1+f(x,y), & x^2+y^2-0.5^2\geq 0,\\
                  f(x,y), & x^2+y^2-0.5^2< 0.
                \end{array}
              \right.
\end{equation}
 It is clear that we can observe some oscillations close to the discontinuities. In order to avoid these non-desirable effects, non-linear methods can be applied when approximating this kind of data. In particular, we reformulate the PUM using Weighted Essentially Non-Oscillatory (WENO) techniques, designed by Liu et al in \cite{Liu} and extensively developed by Shu (see \cite{doi:10.1137/070679065}). The basic idea is to modify the weights of the partition of unity by replacing them with non-linear ones that depend on the data $F$. This way, when a discontinuity crosses a subdomain, the approximation is not considered. It is not difficult to prove that the order of accuracy is conserved in all the domain, except in a band close to the discontinuity.

The paper is divided in the following sections: In Section \ref{rbfinterpolation}, we start by reviewing the RBF interpolation and showing some theorems where the error bounds are estimated, Sec. \ref{pum} is devoted to introduce the PUM and to prove the theoretical results related to the order of accuracy. Afterwards, in Section \ref{weno}, we review WENO method from the PUM perspective, identifying each component of WENO with its corresponding component in PUM. The non-linear method is presented in Sec. \ref{NLPUM} where we explain its theoretical properties and its implementation. Some numerical results are performed in Sec. \ref{expnum} and some conclusions and future work are outlined in the last section.

\section{Radial Basis Function Interpolation}\label{rbfinterpolation}

There is extensive literature on the RBF interpolation method (see e.g. \cite{buhmann,FASSHAUER,wendland}), we will principally employ the notation presented in \cite{cavo,FASSHAUER}. We consider a set of scattered data $X=\{\xx_i\}_{i=1}^N$ in a bounded domain $\Omega\subset \mathbb{R}^n$, with an associated set of data values $F=\{f_i:=f(\xx_i)\}_{i=1}^N$ which are the evaluations of certain unknown function $f:\Omega\to \mathbb{R}$ at the nodes $\xx_i$. The scattered data interpolation problem is the following: to find an operator, $\mathcal{I}(f):\Omega\to \mathbb{R}$, such that:
$$\mathcal{I}(f)(\xx_i)=f(\xx_i),\quad i=1,\hdots,N.$$
Ideally, we aim to minimize the distance between $\mathcal{I}(f)$ and $f$ in a specified norm. First, we define a radial function $\Phi:\Omega\to \mathbb{R}$ when there exists a univariate function $\phi:[0,\infty[\to \mathbb{R}$ such that
$$\Phi(\xx)=\phi(\|\xx\|_2), \quad\forall\,\, \xx\in\Omega,$$
with $\|\cdot\|_2$ the Euclidian norm in $\mathbb{R}^n.$ Now, we denote as $\Phi_i(\xx)=\phi(\|\xx-\xx_i\|_2)$, $i=1,\hdots,N$ and use a radial basis function expansion
to formulate the problem:
\begin{equation}\label{inteporlator}
\mathcal{I}(f)(\xx)=\sum_{i=1}^N c_i \Phi_i(\xx)=\sum_{i=1}^N c_i \phi(\|\xx-\xx_i\|_2),\quad \xx\in\Omega,
\end{equation}
being $c_i$ the coefficients that solve the interpolation system $A\mathbf{c}=\mathbf{f}$ with
$$
A=\begin{bmatrix}
\phi(\| \mathbf{x}_1 - \mathbf{x}_1 \|_2) & \phi(\| \mathbf{x}_1 - \mathbf{x}_2 \|_2) & \cdots & \phi(\| \mathbf{x}_1 - \mathbf{x}_N \|_2) \\
\phi(\| \mathbf{x}_2 - \mathbf{x}_1 \|_2) & \phi(\| \mathbf{x}_2 - \mathbf{x}_2 \|_2) & \cdots & \phi(\| \mathbf{x}_2 - \mathbf{x}_N \|_2) \\
\vdots & \vdots & \ddots & \vdots \\
\phi(\| \mathbf{x}_N - \mathbf{x}_1 \|_2) & \phi(\| \mathbf{x}_N - \mathbf{x}_2 \|_2) & \cdots & \phi(\| \mathbf{x}_N - \mathbf{x}_N \|_2)
\end{bmatrix}, \,\,
\mathbf{c}=\begin{bmatrix}
c_1 \\
c_2 \\
\vdots \\
c_N
\end{bmatrix}, \,\, \mathbf{f}=
\begin{bmatrix}
f_1 \\
f_2 \\
\vdots \\
f_N
\end{bmatrix}.
$$
It is clear that $A$ is symmetric and some conditions over $\phi$ has to be imposed to get a matrix $A$ definite positive and, therefore, a unique solution of the
system. Next definition is useful for this purpose (see \cite{FASSHAUER}).

\begin{definition}[\cite{FASSHAUER}]
A real continuous function $\Phi:\mathbb{R}^n\to\mathbb{R}$ is called positive definite on $\mathbb{R}^n$ if
\begin{equation}\label{formacuadratica}
\sum_{j=1}^N \sum_{i=1}^N c_j c_i \Phi_i(\mathbf{x}_j) \geq 0
\end{equation}
for any $N$ pairwise different points $\{x_i\}_{i=1}^N\subset \mathbb{R}^n$ and $\mathbf{c}=[c_1,\hdots,c_N]\in\mathbb{R}^N$. It is strictly positive definite
on $\mathbb{R}^n$ if the form \eqref{formacuadratica} is zero only for $\mathbf{c}=\boldsymbol{0}$.
\end{definition}

In the literature we can find many examples of functions $\phi$ which produce strictly positive definite functions, such as
Gaussian function $\phi_G(\|\xx\|_2)=e^{-\|\xx\|_2^2}$, that is also $\mathcal{C}^\infty$. Another example are Wendland functions, defined as
\begin{equation}\label{wendlandfunctions}
\begin{split}
&\phi_{W_0}(\|\xx\|_2)=(1-\|\xx\|_2)^2_+,\\
&\phi_{W_2}(\|\xx\|_2)=(1-\|\xx\|_2)^4_+(4\|\xx\|_2+1),\\
&\phi_{W_4}(\|\xx\|_2)=(1-\|\xx\|_2)^6_+(35\|\xx\|^2_2+18\|\xx\|_2+3),
\end{split}
\end{equation}
based on the cutoff function
\[(x)_+=\begin{cases}
x & \text{for}\,\, x\geq 0,\\
0 & \text{for}\,\, x<0 ,
\end{cases}
\]
which are $\mathcal{C}^0$, $\mathcal{C}^2$ and $\mathcal{C}^4$ respectively. Mat\'ern functions are another example, that are defined as
\begin{equation}\label{maternfunctions}
\begin{split}
&\phi_{M_0}(\|\xx\|_2)=e^{-\|\xx\|_2},\\
&\phi_{M_2}(\|\xx\|_2)=(1 + \|\xx\|_2) e^{-\|\xx\|_2},\\
&\phi_{M_4}(\|\xx\|_2)=\left(3 + 3\|\xx\|_2 + \|\xx\|_2^2\right) e^{-\|\xx\|_2},
\end{split}
\end{equation}
satisfying that $\phi_{M_k}\in \mathcal{C}^k$, $k=0,2,4$. Other classes of functions, such as
 generalized inverse multiquadrics or Poisson radial functions, also produce strictly positive definite functions, see \cite{FASSHAUER} for more details.

The first characteristic of this interpolant,  Eq. \eqref{inteporlator}, is the smoothness, as it is a linear combination of smooth functions, $\Phi_i$. Consequently, the continuity of the operator only depends on the continuity of the function $\phi(\|\cdot\|_2)$.

Secondly, the error estimates and the order of accuracy for the interpolation process have been studied extensively (see \cite{buhmann,FASSHAUER,wendland} and the references cited therein for further details). To introduce the results, we need to present three relevant concepts: First of all, we define the fill distance, which measures the extent to which the data in the set \( X \) adequately cover the domain \( \Omega \), as
$$h=h_{X,\Omega}=\sup_{x\in \Omega}\min_{x_j\in X}\|x-x_j\|_2.$$
Following this, we refer to the space (see \cite{cavo16,FASSHAUER})
$$\mathcal{C}^{2k}_\nu(\mathbb{R}^n)=\{f\in \mathcal{C}^{2k}:D^{\boldsymbol{\beta}}f(\xx)=O(\|\xx\|_2^\nu),\text{ as }\|\xx\|_2\to 0\,\,\text{ if } |\boldsymbol{\beta}|=2k\},$$
and finally, we state that a region $\Omega\subset \mathbb{R}^n$ satisfies \textbf{an interior cone condition}, \cite{FASSHAUER,wendland}, if there exists an angle $\theta\in (0,\pi/2)$ and a radius
$r>0$ such that for every $\xx\in\Omega$ there exists a unit vector $\boldsymbol{\xi}(\xx)$ such that the cone
$$C=\{\xx+\lambda\mathbf{y}:\mathbf{y}\in\mathbb{R}^n,\,\, \|\mathbf{y}\|_2=1,\,\, \mathbf{y}^T\boldsymbol{\xi}(\xx)\geq \cos \theta, \,\,\lambda\in[0,r]\}.$$
The following lemma presented by Wendland in \cite{wendland} indicates that a ball always satisfies an interior cone condition. This result will be relevant in our numerical experiments when we construct a particular partition of unity method.
\begin{lemma}
Every ball with radius \(\delta > 0\) satisfies an interior cone condition with radius \(\delta > 0\) and angle \(\theta = \pi / 3\).
\end{lemma}

With these three elements: the fill distance, the space $\mathcal{C}^{2k}_\nu(\mathbb{R}^n)$ and the interior cone condition,  we review the following result.
\begin{theorem}\label{thrbf}{\cite{FASSHAUER}}
Suppose $\Omega\subset \mathbb{R}^n$ is bounded and satisfies an interior cone condition. Suppose $\Phi\in \mathcal{C}^{2k}_\nu(\mathbb{R}^n)$ strictly positive definite. Then, the error between $f\in \mathcal{N}_{\Phi}(\Omega):=\text{\emph{span}}\{\Phi(\cdot-\xx),\xx\in X\}$ and its interpolant $\mathcal{I}(f)$ defined in Eq. \eqref{inteporlator} can be bounded by
$$|f(\xx)-\mathcal{I}(f)(\xx)|\leq Ch^{k+\frac{\nu}{2}}\|f\|_{\mathcal{N}_\Phi(\Omega)}$$
for all $\xx\in \Omega$, $h$ is sufficiently small and where $C$ is a constant independent of $\xx$, $f$ and $\Phi$.
\end{theorem}

In this work, we focus on studying the differences between the function and the interpolator. Error estimates for the derivatives are proved in \cite{FASSHAUER,wendland}. To conclude this review section, we introduce
the following results for Mat\'ern function, Eq. \eqref{maternfunctions}, 
as they will be used in the numerical experiments. For this purpose, we need to define the Sobolev spaces $W^m_2(\mathbb{R}^n)$ introduced in \cite{FASSHAUER} as
$$W^{m}_2(\Omega)=\{f\in L_2(\Omega)\cap\mathcal{C}(\Omega):D^{\boldsymbol{\beta}}f \in L_2(\Omega)\,\,\text{ for all }\,\,|\boldsymbol{\beta}|\leq m,\,\,\boldsymbol{\beta}\in\mathbb{N}^n\}.$$

\begin{proposition}{\cite{FASSHAUER}}
For the strictly positive Mat\'ern functions $\phi_{M_k}$, $k=0,2,4$, Eq. \eqref{maternfunctions}, is satisfied that
$$|f(\xx)-\mathcal{I}_{M_k}(f)(\xx)|\leq Ch^{\frac{k+1}{2}}\|f\|_{W^{1+k/2}_2(\Omega)}$$
for all $\xx\in \Omega$, $h$ is sufficiently small, and $f\in W^{1+k/2}_2(\Omega)$.
\label{propmatern}
\end{proposition}

Next subsection is devoted to introduce the partition of unity method for radial basis functions, which is a well-known approximation method
that is satisfactorily used when the number of data points is very large.

\section{Partition of Unity Approximation}\label{pum}

The partition of unity method (PUM) is proposed to facilitate efficient computation using meshfree approximation methods. The approach is straightforward: it involves decomposing the large problem into smaller subproblems while preserving the order of accuracy. Accordingly, using the same notation as above, let \( \Omega \subset \mathbb{R}^n \) represent our open and bounded domain, which we partition into $M$ subdomains $\Omega_j$ such that \cite{cavo16,wendland2}:
\begin{itemize}
\item The domain is contained in the union of the subdomains, i.e. $$\Omega\subseteq \bigcup_{j=1}^M \Omega_j.$$
\item For every $\xx\in\Omega$ the number of subdomains $\Omega_j$ with $\xx\in\Omega_j$ is bounded by a global constant $K$.
\item There exists a constant \( C_r > 0\) and an angle \( \theta \in (0, \pi/2) \) such that every patch \( \Omega_j \cap \Omega \) satisfies an interior cone condition with angle \( \theta \) and radius \( r = C_r h_{X,\Omega} \).
\item The local fill distance $h_{X_j,\Omega_j}$ with $X_j=\Omega_j\cup X$ are uniformly bounded by the global fill distance $h_{X,\Omega}$.
\end{itemize}
This is called a regular covering for $(\Omega,X)$, \cite{wendland}. For this covering, we choose a partition of unity, i.e. a family of compactly supported, non-negative, continuous functions, $w_j:\Omega_j\to \mathbb{R}$ with
$$\sum_{j=1}^M w_j(\xx)=1\,\,\forall\,\,\xx\in\Omega,\quad \text{supp}(w_j)\subseteq \Omega_j,$$
and define the global interpolator as:
\begin{equation}\label{puinterpolator}
\mathcal{I}_{\text{PUM}}(f)(\xx)=\sum_{j=1}^M w_j(\xx)\mathcal{I}_j(f)(\xx),
\end{equation}
where $\mathcal{I}_j(f):\Omega_j\to \mathbb{R}$ is a local RBF interpolant on each sudomain $\Omega_j$.

In \cite{FASSHAUER} some technical conditions over the partition of unity functions are imposed to establish error bounds. We assume that the
partition of unity functions is $k$-stable. This is, each $w_j\in \mathcal{C}^k(\mathbb{R}^n)$ and for every multi-index $\boldsymbol{\alpha}$ with
$|\boldsymbol{\alpha}|\leq k$ the following inequality is satisfied:
$$\|D^{\boldsymbol{\alpha}} w_j\|_{\infty,\Omega_j}\leq \frac{C_{\boldsymbol{\alpha}}}{\delta_j^{|\boldsymbol{\alpha}|}},$$
where $C_{\boldsymbol{\alpha}}$ is some positive constant and $\delta_j$ is the diameter of $\Omega_j$, i.e. $\delta_j=\sup_{x,y\in\Omega_j}\|x-y\|_2$. With these restrictions over the covering and the partition of unity function, we show the following result, \cite{FASSHAUER}.
\begin{theorem}\label{thpum}
Suppose \(\Omega \subseteq \mathbb{R}^n\) is open and bounded, and let \(X = \{\xx_i\}_{i=1}^N \subseteq \Omega\). Let \(\Phi \in \mathcal{C}_\nu^{2k}(\mathbb{R}^n)\) be strictly positive definite. Let \(\{\Omega_j\}_{j=1}^M\) be a regular covering for \((\Omega, {X})\) and let \(\{w_j\}_{j=1}^M\) be \(k\)-stable for \(\{\Omega_j\}_{j=1}^M\). Then the error between \(f \in \mathcal{N}_\Phi(\Omega)\) and its partition of unity interpolant, Eq. \eqref{puinterpolator}, can be bounded by
\[
|f(\xx) - \mathcal{I}_{\text{\emph{PUM}}}(f)(\xx)| \leq C h^{k + \frac{\nu}{2}} \|f\|_{\mathcal{N}_\Phi(\Omega)},
\]
for all \(\xx \in \Omega\).
\end{theorem}
It is clear that the order of accuracy does not vary using PUM or RBF method, as we can see in the theorems \ref{thrbf} and \ref{thpum}. A simple way, suggested in \cite{cavo16,FASSHAUER}, to construct PUM is to consider Shepard's approximation, i.e.
\begin{equation}\label{wjs}
w_j(\xx)=\frac{\varphi_j(\xx)}{\sum_{k=1}^M\varphi_k(\xx)},\quad j=1,\hdots,M,
\end{equation}
where $\varphi_j$ are compactly supported, with their support contained in $\Omega_j$, and can be chosen, for example, as Wendland functions, Eq. \eqref{wendlandfunctions}.

\section{Reviewing Weighted Essentially Non-Oscillatory algorithm as a PUM}\label{weno}

The WENO scheme has gained widespread recognition as an effective method for approximating data values, especially in contexts where discontinuities are prevalent (see e.g. \cite{doi:10.1137/070679065}). Initially introduced for approximatig the solution of hyperbolic partial differential equations, the WENO approach is structured to achieve high-order accuracy while mitigating the unwanted oscillations that typically occur near discontinuities \cite{JiangShu}. This is accomplished through the adaptive assignment of weights that prioritize smoother portions of the data and diminish the contribution of stencils that intersect discontinuities.

In what follows, we review the WENO-$2r$ method in the point-values context (see \cite{ARSY20,ABM})  relating each of its component to those corresponding to the PUM. The analogy is done for the algorithms in one variable, but it can be generalized to several variables \cite{muletruizshuyanez}. In this case, let $\{x_i\}^{J}_{i=0}$ be a uniform partition of the interval
$[a, b]$ in $J$ subintervals,
$x_i=a+i\cdot h, \quad h=\frac{b-a}{J},$ and let $\left\{f_i\right\}_{i=0}^{J}$ be the point-values discretization of an unknown function $f$ at the nodes $x_i$, to interpolate at $x\in (x_{i-1}, x_i)$. WENO-$2r$ algorithm uses the stencil $X=\{x_{i-r},\cdots, x_{i+r-1}\}$, that is composed of $N=2r$ nodes, we denote it by $\tilde{\mathcal{I}}(f)$. We can divide the stencil $X$ in $r$ substencils of $r+1$ nodes, $\mathcal{S}_j=\{x_{i-r+j},\hdots,x_{i+j}\}$, $j=0,\hdots,r-1$ and prove that there exist $\tilde{w}_j:[x_{i-1}, x_i]\to \mathbb{R}$, $j=0,\hdots,r-1$, \cite{muletruizshuyanez}, called optimal weights, such that:
$$\tilde{\mathcal{I}}(f)(x)=\sum_{j=0}^{r-1} \tilde{w}_j(x) \tilde{\mathcal{I}}_j(f)(x),$$
with $\tilde{w}_j(x)>0$ and $\sum_{j=0}^{r-1}\tilde{w}_j(x)=1$. This is a PU-like method, with $\Omega_j=[x_{i-r+j},x_{i+j}]$ and $\tilde{\mathcal{I}}_j(f)$ the Lagrange interpolatory polynomial using the nodes $\mathcal{S}_j$. Now, the WENO method consists in replacing $\tilde{w}_j(x)$ by non-linear weights $\tilde{\mathcal{W}}_j(x)$, such that
\begin{equation}
\mathcal{I}_{\text{WENO}}(f)(x) = \sum_{j=0}^{r-1} \tilde{\mathcal{W}}_j(x)\tilde{\mathcal{I}}_j(f)(x),
\end{equation}
where the non-linear weights are defined as:
\begin{equation}\label{Wjs}
\tilde{\mathcal{W}}_j(x) = \frac{\tilde{\alpha}_j(x)}{\sum_{k=0}^{r-1} \tilde{\alpha}_k(x)},
\,\,\text{with}\,\,
\tilde{\alpha}_j(x) = \frac{\tilde{w}_j(x)}{(\epsilon + \tilde{I}_j)^t},
\end{equation}
being $\tilde{I}_j$ smoothness indicators which indicates if a discontinuity crosses $\Omega_j$, $\epsilon$ is a small positive parameter to avoid division by zero, and $t$ is typically a constant to enhance accuracy in smooth regions. The idea is to make a non-linear version of the interpolator $\mathcal{I}_{\text{PUM}}$, Eq. \eqref{puinterpolator}, (NL-PUM) using the ideas of WENO method. The key is to define an adequate smoothness indicator for our framework.

\section{Non-linear partition of unity method}\label{NLPUM}

As we mentioned in the previous section, Sec. \ref{weno},  we construct a non-linear PUM starting with the interpolator $\mathcal{I}_{\text{PUM}}(f)$, Eq. \eqref{puinterpolator}:
\begin{equation*}
\mathcal{I}_{\text{PUM}}(f)(\xx)=\sum_{j=1}^M w_j(\xx)\mathcal{I}_j(f)(\xx),
\end{equation*}
replacing the weight functions $w_j$ described in Eq. \eqref{wjs} by non-linear ones in the same way as in Eq. \eqref{Wjs}, i.e.
\begin{equation*}
\mathcal{W}_j(\xx) = \frac{\alpha_j(\xx)}{\sum_{k=1}^{M} \alpha_k(\xx)},
\,\,\text{with}\,\,
\alpha_j(\xx) = \gamma_j \varphi_j(\xx)=\frac{\varphi_j(\xx)}{(\epsilon + I_j)^t},
\end{equation*}
where $0<\gamma_j=(\epsilon + I_j)^{-t}$ is a positive constant, and $\varphi_j$ are compactly supported functions with support on $\Omega_j$. Note that the new functions $\alpha_j$ are compactly supported and, if we choose $\varphi_j$ as Wendland functions, then $\alpha_j$ are Wendland functions too. The constant $\epsilon=10^{-14}$ is set to avoid non-zero denominators. The parameter $t$ is relevant to correctly detect the discontinuities and we choose $t=6$. Finally, $I_j$, $j=1,\hdots, M$ are the smoothness indicators. To construct them, we consider
$$\mathcal{S}_j=\Omega_j\cap X =\{\xx_{j_i}\in \Omega_j\cap X :i=1,\hdots,N_j\},$$
with $N_j>3$, and the linear least square problem
$$p_j=\argmin_{p\in\Pi_1(\mathbb{R}^2)} \sum_{i=1}^{N_j} (f(\xx_{j_i})-p(\xx_{j_i}))^2,$$
where $\Pi_1(\mathbb{R}^2)$ is the set of polynomials of degree less than or equal to one, and define
\begin{equation}\label{smoothindicators}
I_j=\frac{1}{N_j}\sum_{i=1}^{N_j}|f(\xx_{j_i})-p_j(\xx_{j_i})|,
\end{equation}
then $I_j$ satisfies the conditions:
\begin{itemize}
\item The order of the smoothness indicator that are free of discontinuities is $h^2$, i.e.
$$I_j=O(h^2) \,\, \text{if}\,\, f \,\, \text{is smooth in } \,\, \Omega_j.$$
\item When a discontinuity crosses $\Omega_j$ then:
$$I_j \nrightarrow 0 \,\, \text{as}\,\, h\to 0.$$
\end{itemize}
Therefore, the new non-linear operator will be:
\begin{equation}\label{nlpum}
\mathcal{I}_{\text{NL-PUM}}(f)(\xx)=\sum_{j=1}^M \mathcal{W}_j(\xx)\mathcal{I}_j(f)(\xx).
\end{equation}
The key to this new method is that it gives more weight to the data placed in smooth regions of the data, neglecting data located near the discontinuity, as this data is multiplied by a weight of order $O(h^{2t})$. In the next subsection we will analyze the properties of the interpolator.

\subsection{Properties of the new NL-PUM}
From the new interpolator, Eq. \eqref{nlpum} is defined as PUM multiplying the compactly supported functions by constants (that are dependent on the data $F=\{f_i\}_{i=1}^N$). The smoothness of the new operator is similar to PUM, int he sense that its smoothness is solely determined by the smoothness of the functions $\varphi(\|\cdot\|_2)$ and $\phi(\|\cdot\|_2)$.

Secondly, we determine the order of accuracy when we apply $\mathcal{I}_{\text{NL-PUM}}(f)$ to data that comes from a function which is piecewise continuous throughout its subdomains but not continuous, i.e. we suppose an open and bounded domain $\Omega\subset\mathbb{R}^n$ and consider a curve $\zeta:\mathbb{R}^n\to \mathbb{R}$, $\zeta\in \mathcal{C}^1$ and the sets
$$\Gamma=\{\xx\in \Omega:\zeta(\xx)=0\},\quad \Omega^+=\{\xx\in \Omega:\zeta(\xx)\geq 0\},\quad \Omega^-=\Omega\setminus \Omega^+$$
that is to say that $\Omega$ is divided into two subsets, $\Omega=\Omega^+\cup \Omega^-$.  Let us assume that our data originates from a function of the type:
$$
f(\xx)=\begin{cases}
f_+(\xx), & \xx \in \Omega_+,\\
f_-(\xx), & \xx \in \Omega_-,\\
\end{cases}
$$
with $f_\pm\in \mathcal{C}^k(\overline{\Omega_\pm})$, $k\in \mathbb{N}$, but
$$\lim_{\xx\to \xx_0 \in \Gamma} f_-(\xx)\neq f_+(\xx_0),$$
then our scattered data $X=\{\xx_i\}_{i=1}^N$ are divided in two $X_\pm=X\cap \Omega_\pm$. We select our partition $\{\Omega_j\}_{j=1}^M$ as in Section \ref{pum} and define patches that are contaminated by a discontinuity as follows.

\begin{definition}[Contaminated patches]
Using the notation presented above, $\Omega_k\in\{\Omega_j\}_{j=1}^M$ is a contaminated patch if there exist
$\Omega_k\cap X_{+}\neq \emptyset$ and $\Omega_k\cap X_{-}\neq \emptyset$ .\label{contapatches}
\end{definition}

\begin{figure}[hbtp!]
\begin{center}
		\begin{tabular}{cc}
			\includegraphics[width=6cm]{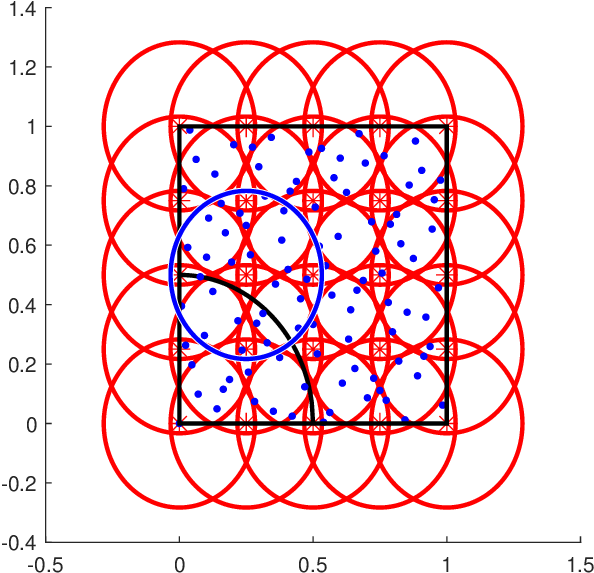} &
			\includegraphics[width=6cm]{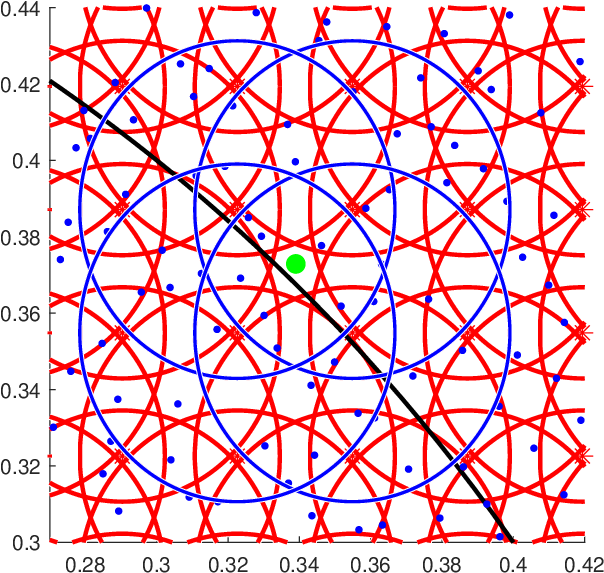}\\
		\end{tabular}
\end{center}
\caption{Black square: frontier of $\Omega$, red circumferences: frontiers of the patches, black curve: discontinuity, blue circumferences: contaminated patches, blue points: data, red stars: centers of the patches, green point (right): discontinuity point.}
	\label{figureejemplot1}
	\end{figure}

We determine if a subset $\Omega_k \in \{\Omega_j\}_{j=1}^M$ is a contaminated patch if $I_k>h=h_{X,\Omega}$. Also, the points close to the discontinuity will be
marked using the following definition.
\begin{definition}[Discontinuity points]
Using the notation above presented, let $\xx \in \Omega$ be a point of the domain, $\varepsilon>0$ be a fixed threshold  and
\begin{equation*}
\begin{split}
&S(\xx)=\{j\in \{1,\hdots,M\}:\xx\in \text{\emph{supp}}(\Omega_j),\, \varphi_j(\xx)>\varepsilon\},\\
&\tilde{S}(\xx)=\{j\in \{1,\hdots,M\}:\xx\in \text{\emph{supp}}(\Omega_j), \,\Omega_k\cap X_{+}\neq \emptyset,\,\,\Omega_k\cap X_{-}\neq \emptyset,\, \varphi_j(\xx)>\varepsilon\},
\end{split}
\end{equation*}
then $\xx$ is a discontinuity point if
$$S(\xx)=\tilde{S}(\xx).$$\label{discpoint}
\end{definition}
This definition means that a point is a discontinuity point only if all the patches used to interpolate at that point are contaminated, Def. \ref{contapatches}.
Last definition, Def. \ref{discpoint}, is quite relevant to design the new algorithm. This is due to the fact that, at the discontinuity points, the strategy to interpolate has to be changed due to the oscillations produced by the RBF method close to the discontinuity curve. Hence, the most problematic case is the one in which all the patches containing the point where the interpolation is desired, $\xx$, are contaminated, except for one whose weight, $\varphi_j(\xx)$, is close to zero. For this reason, we impose the threshold condition in the previous definition,  Def. \ref{discpoint}. In Fig. \ref{figureejemplot1} we can see a picture with the representation of the two previous concepts, Defs. \ref{contapatches} and \ref{discpoint}. We represent the data points using blue dots, the patch boundaries with red circles, the discontinuity curve and the domain boundary with black lines, and the patch centers with red stars. In the figure on the left, the contaminated patch is indicated with a blue circle. In the figure on the right, the point at which the function is to be approximated is marked in green, while the patches used for interpolation are highlighted in blue. It can be observed that all the patches are contaminated, indicating that the green point is a discontinuity point.

With these definitions, the following result is straightforward based on Th. \ref{thpum}, where an error bound is deduced dependent on the smoothness of
the function $\Phi$ and its native space $\mathcal{N}_\Phi$.
\begin{theorem}\label{thnlpum}
Suppose \(\Omega \subseteq \mathbb{R}^n\) is open and bounded, and let \(X = \{\xx_i\}_{i=1}^N \subseteq \Omega\). Let \(\Phi \in \mathcal{C}_\nu^{2k}(\mathbb{R}^n)\) be strictly positive definite. Let \(\{\Omega_j\}_{j=1}^M\) be a regular covering for \((\Omega, {X})\) and let \(\{w_j\}_{j=1}^M\) be \(k\)-stable for \(\{\Omega_j\}_{j=1}^M\). Then, if $\xx\in \Omega_+$, it is not a discontinuity point, Def. \eqref{discpoint}, and \(f_\pm \in \mathcal{N}_\Phi(\Omega_\pm)\), the error between \(f\) and its non-linear partition of unity interpolant, Eq. \eqref{nlpum}, with $t\geq k/2+\nu/4$ can be bounded by
\[
|f(\xx) - \mathcal{I}_{\text{\emph{NL-PUM}}}(f)(\xx)| \leq C h^{k + \frac{\nu}{2}} \|f\|_{\mathcal{N}_\Phi(\Omega_+)}.
\]
Analogously, if $\xx \in \Omega_-$ and it is not a discontinuity point.
\end{theorem}

In other words, if $\xx\in \Omega$ and there exists, at least, a patch $\Omega_k$ free of discontinuities, then the order of accuracy is conserved. Consequently, the discontinuity points have to be treated in a special manner. Also, if we suppose (without loss of generality) that
$$\lim_{\xx\to\xx_0\in \Gamma} f_+(\xx) < \lim_{\xx\to\xx_0\in \Gamma} f_-(\xx), $$
we can not ensure that if $\tilde{\xx}\in \Omega$ is a discontinuity point then
$$\lim_{\xx\to\xx_0\in \Gamma} f_+(\xx) \leq  \mathcal{I}_{\text{{NL-PUM}}}(f)(\tilde{\xx}) \leq   \lim_{\xx\to\xx_0\in \Gamma} f_-(\xx). $$
For this reason, in order to avoid some oscillations of the interpolator  $\mathcal{I}_{\text{{NL-PUM}}}(f)$ close to the discontinuities, we modify our algorithm and explain it in the next subsection.

\subsection{On the implementation of the method. Avoiding oscillations close to the discontinuities}

We compute our method using the implementation explained by Cavoretto et al. in \cite{cavo}, based on the partitioning data procedures proposed in \cite{cavo15}. Let $\{\hat{\xx}_k\}_{k=1}^{L}\subset \Omega$ be a sample where we want to approximate the function $f$,
we start with the covering of the domain $\Omega$, consider some centers $P=\{\bar{\xx}_j\}_{j=1}^M\subset\Omega$, fix a radius $\delta$ and construct the subdomains $\Omega_j =B(\bar{\xx}_j,\delta)$ where $B(\bar{\xx}_j,\delta)$ is a ball of center $\bar{\xx}_j$ and radius $\delta$. It is important to remark, as noted in \cite{cavo,FASSHAUER}, that if a nearby node distribution is chosen, $M$ is an adequate number of PU subdomains if
$$\frac{N}{M} \approx 2^n,$$
then, the radius has to satisfy the condition, \cite{cavo}:
\begin{equation}\label{deltacondition}
\delta\geq \frac{1}{M^\frac1n}.
\end{equation}
At this point, we estimate the smoothness indicators $I_j$, Eq. \eqref{smoothindicators} and mark the contaminated patches, Def. \ref{contapatches}. The interpolation data located in each subdomain $j$ are found using the partitioning data structures explained in \cite{cavo}. After this, we estimate the non-linear weights in each point $\hat{\xx}_k$  using the smoothness indicators calculated, and label the discontinuity points, we called them $\tilde{\xx}_j$. We apply the classical Shepard's method only for these points, with the same compactly supported function used in PUM. In summary, we describe the algorithm in two passes:
\begin{enumerate}
\item We use the software designed in \cite{cavo} and mark the discontinuity points, $\tilde{\xx}$.
\item We apply the Shepard's method to these discontinuity points.
\end{enumerate}
Due to the partitioning data structures introduced in \cite{cavo}, which allow us to organize data very efficiently across different subdomains, the modifications inserted to make the algorithm non-linear are not too costly in terms of computational efficiency.

\section{Numerical experiments}\label{expnum}

In this section, we perform some tests in order to check two principal properties of the new non-linear algorithm, NLPUM, in comparison with the classical one, PUM. We will conduct the experiments using a uniform grid and Halton's points as data. Firstly, in Subsection \ref{orderofaccuracy}, we employ our method with a continuous function, in particular with Franke's function, Eq. \eqref{frankesfunction}, to verify the order of accuracy. Next subsection is devoted to analyze the elimination of the non-desired oscillations close to the discontinuities, we come back to the example presented in the introduction, Sec. \ref{introduction}, and show the powerful capabilities of the new algorithm for this type of functions. We take the same parameters as \cite{cavo} in all the experiments, i.e., the number of patches in one direction is $\lfloor \sqrt{N}/2 \rfloor$ and the radius $\delta=\sqrt{2/M}$.

\subsection{Order of accuracy}\label{orderofaccuracy}

We begin by analyzing the order of accuracy using the widely recognized Franke's function, Eq. \eqref{frankesfunction}, using as data points a regular grid defined by \(X_{N=2^l+1}=\{(i/2^l,j/2^l):i,j=0,\hdots,2^l\}\), and a set of \(N=(2^l+1)^2\) Halton's scattered data points, \cite{halton}. We denote the fill distance by \(h_l\), and the errors by \(e^l_k=|f(\hat{\mathbf{x}}_k)-\mathcal{I}^l(\hat{\mathbf{x}}_k)|\), where \(\{\hat{\mathbf{x}}_k\}_{k=1}^{60}\) is a regular grid in $[0,1]^2$.
We define the maximum discrete \(\ell^2\) norms and their respective convergence rates as follows:
\begin{equation*}
\begin{split}
&\text{MAE}_l=\max_{i=1,\hdots,60}e^l_i, \quad \text{RMSE}_l=\left(\frac{1}{60}\sum_{i=1}^{60}(e_i^l)^2\right)^{\frac12}, \\
&r_l^{\infty}=\frac{\log(\text{MAE}_{l-1}/\text{MAE}_{l})}{\log(h_{l-1}/h_{l})}, \quad r_l^2=\frac{\log(\text{RMSE}_{l-1}/\text{RMSE}_{l})}{\log(h_{l-1}/h_{l})}.
\end{split}
\end{equation*}
In the first example, we choose Mat\`ern function, $\phi_{M_2}$ for the RBF problems and the Wendland $\mathcal{C}^2$ function $\phi_{W_2}$, for the partition of unity. We present the results in Table \eqref{tabla:d0h}, and it is clear that they are very similar for  both the linear and non-linear methods. The expected theoretical order of accuracy, Prop. \eqref{propmatern}, is improved for the two algorithms in this example. With respect to the errors, the NLPUM yields slightly better results when the data is placed on an uniform grid.

\begin{table}[hbpt!]
\begin{center}
\begin{tabular}{lrrrrrrrrrrr}
\hline  \multicolumn{5}{c}{Grid points} & &   \multicolumn{5}{c}{Halton's points}\\\hline
& \multicolumn{2}{c}{PUM} & &   \multicolumn{2}{c}{NLPUM}& & \multicolumn{2}{c}{PUM} & &   \multicolumn{2}{c}{NLPUM}\\ \cline{2-3} \cline{5-6} \cline{8-9} \cline{11-12}
$l$  & $\text{MAE}_l$ & $r_l^{\infty}$ & &$\text{MAE}_l$ & $r_l^{\infty}$         & &  $\text{MAE}_l$ & $r_l^{\infty}$ & &$\text{MAE}_l$ & $r_l^{\infty}$           \\\hline
$4$  &           7.5944e-03 & -            & &       2.8436e-01 & -               & &         1.4157e-02 & -            & &        2.6937e-02 & -                 \\
$5$  &           1.6541e-03 &       2.1989 & &       8.6890e-02 &          1.7105 & &         2.6921e-03 &      2.1831 & &         2.9567e-03&             2.9058   \\
$6$  &           3.7494e-04 &       2.1413 & &       8.8016e-05 &          9.9472 & &         9.1058e-04 &      2.3481 & &         6.6469e-04&             3.2331   \\
$7$  &           1.1278e-04 &       1.7332 & &       2.6376e-05 &          1.7385 & &         2.0617e-04 &      1.9743 & &         3.0466e-04&             1.0369   \\
$8$  &           2.8906e-05 &       1.9640 & &       7.0429e-06 &          1.9050 & &         6.0296e-05 &      1.4750 & &         7.9830e-05&             1.6068   \\
\hline
\hline
& \multicolumn{2}{c}{PUM} & &   \multicolumn{2}{c}{NLPUM}& & \multicolumn{2}{c}{PUM} & &   \multicolumn{2}{c}{NLPUM}\\ \cline{2-3} \cline{5-6} \cline{8-9} \cline{11-12}
$l$   &$\text{RMSE}_l$ & $r_l^2$ & &$\text{RMSE}_l$ & $r_l^2$ &  & $\text{RMSE}_l$ & $r_l^2$ & &$\text{RMSE}_l$ & $r_l^2$    \\\hline
$4$  &           1.3800e-03 & -         & &      9.6360e-02& -            &&      6.8732e-03 & -         & &      9.6360e-02& -             \\
$5$  &           2.2904e-04 &     2.5910& &      4.7798e-03&       4.3334 &&      3.3489e-04 &      2.2581& &      4.7798e-03&      3.9740  \\
$6$  &           5.8939e-05 &     1.9583& &      1.5566e-05&       8.2624 &&      3.8497e-05 &      2.7589& &      1.5566e-05&      4.6860  \\
$7$  &           1.8977e-05 &     1.6350& &      1.9141e-06&       3.0236 &&      5.9659e-06 &      1.8766& &      1.9141e-06&      2.4782  \\
$8$  &           4.0380e-06 &     2.2325& &      2.5281e-07&       2.9206 &&      7.4670e-07 &      1.7219& &      2.5281e-07&      2.4933  \\
\hline
\end{tabular}
\end{center}
\caption{Errors and rates using PUM and NLPUM for Franke's test function evaluated at grid points, using $\phi_{M_2}$ for RBF problems, and $\phi_{W_2}$ for partition of unity.}\label{tabla:d0h}
\end{table}
To finish this subsection, we repeat the experiment using the functions $\phi_{M_4}$ for the RBF method and $\phi_{W_4}$ for the partition of unity. The results are
shown in Table \ref{tabla:d1h}. Again, we can see that the errors and orders of accuracy are very similar.

\begin{table}[hbpt!]
\begin{center}
\begin{tabular}{lrrrrrrrrrrr}
\hline  \multicolumn{5}{c}{Grid points} & &   \multicolumn{5}{c}{Halton's points}\\\hline
& \multicolumn{2}{c}{PUM} & &   \multicolumn{2}{c}{NLPUM}& & \multicolumn{2}{c}{PUM} & &   \multicolumn{2}{c}{NLPUM}\\ \cline{2-3} \cline{5-6} \cline{8-9} \cline{11-12}
$l$  & $\text{MAE}_l$ & $r_l^{\infty}$ & &$\text{MAE}_l$ & $r_l^{\infty}$         & &  $\text{MAE}_l$ & $r_l^{\infty}$ & &$\text{MAE}_l$ & $r_l^{\infty}$           \\\hline
$4$  &             6.4397e-03 & -              & &        2.9103e-01 & -              & &            8.1704e-03 & -           & &          8.7701e-03& -                 \\
$5$  &             7.6400e-04 &         3.0754 & &        9.0838e-02 &         1.6798 & &            8.9459e-04 &      2.9091 & &           1.1375e-03&              2.6863   \\
$6$  &             7.3358e-05 &         3.3805 & &        2.0297e-05 &        12.1278 & &            3.0378e-04 &      2.3397 & &           1.6615e-04&              4.1672   \\
$7$  &             9.0336e-06 &         3.0216 & &        1.6754e-06 &         3.5986 & &            2.2676e-05 &      3.4490 & &           4.5792e-05&              1.7129   \\
$8$  &             1.5305e-06 &         2.5613 & &        2.5780e-07 &         2.7002 & &            2.3915e-06 &      2.6987 & &           1.0517e-05&              1.7650   \\
\hline
\hline
& \multicolumn{2}{c}{PUM} & &   \multicolumn{2}{c}{NLPUM}& & \multicolumn{2}{c}{PUM} & &   \multicolumn{2}{c}{NLPUM}\\ \cline{2-3} \cline{5-6} \cline{8-9} \cline{11-12}
$l$   &$\text{RMSE}_l$ & $r_l^2$ & &$\text{RMSE}_l$ & $r_l^2$ &  & $\text{RMSE}_l$ & $r_l^2$ & &$\text{RMSE}_l$ & $r_l^2$    \\\hline
$4$  &           1.0412e-03 & -           & &     1.0168e-01& -             &&       9.6168e-04 & -         & &        2.6447e-03& -             \\
$5$  &           1.0970e-04 &       3.2465& &     4.7209e-03&        4.4289 &&       1.0395e-04 &        2.9260& &      1.2230e-04&        4.0427  \\
$6$  &           1.2027e-05 &       3.1893& &     2.7164e-06&       10.7631 &&       1.3923e-05 &        4.3550& &      6.5149e-06&        6.3522  \\
$7$  &           1.8291e-06 &       2.7170& &     1.5998e-07&        4.0858 &&       1.7160e-06 &        2.7826& &      6.6212e-07&        3.0389  \\
$8$  &           1.9361e-07 &       3.2400& &     1.0651e-08&        3.9088 &&       2.0360e-07 &        2.5573& &      5.6579e-08&        2.9512  \\
\hline
\end{tabular}
\end{center}
\caption{Errors and rates using PUM and NLPUM for Franke's test function evaluated at grid points, using $\phi_{M_4}$ for RBF problems, and $\phi_{W_4}$ for partition of unity.}\label{tabla:d1h}
\end{table}

\subsection{Avoiding discontinuities}

In this subsection, we conduct a series of experiments to analyze the behavior of the new interpolator near discontinuities. In all the experiments the initial number of data points is $N=65^2$.

We begin with a set of experiments in which the original domain $[0,1]^2$ is divided into two subdomains. We assume that the data originates from Franke's function $f$, as defined in Eq. \eqref{frankesfunction}, in one subdomain, and from $1 + f$ in the other, i.e., we define the functions:
\begin{equation*}
f_q(x,y)=\left\{
                \begin{array}{ll}
                  1+f(x,y), & (x,y)\in \Omega^+, \\
                  f(x,y), & (x,y)\in \Omega^-,
                \end{array}
              \right.
\end{equation*}
where $\Gamma_q=\{(x,y)\in [0,1]^2:\zeta_q(x,y)=0\}$, $\Omega^+=\{\xx\in [0,1]^2:\zeta_q(\xx)\geq 0\}$, $\Omega^-=[0,1]^2\setminus \Omega^+$. In this series of experiments we employ $\phi_{M_2}$ as function for RBF  problems and $\varphi=\phi_{W_2}$ for the partition of unity.

\begin{figure}[hbpt!]
\begin{center}
		\begin{tabular}{cc}
PUM & NL-PUM\\
			 \hspace{-1cm}	\includegraphics[width=8cm]{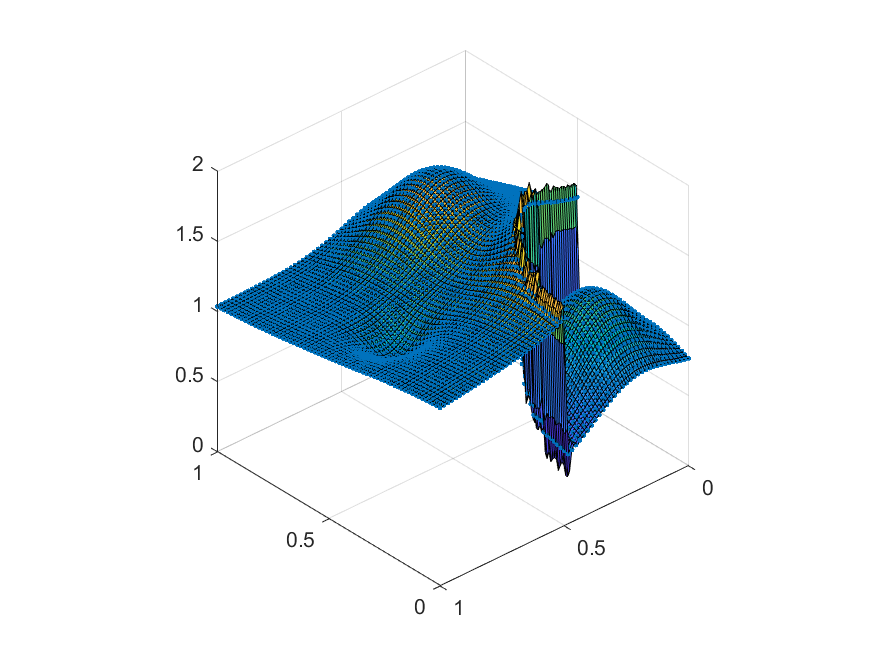} &
		 \hspace{-1cm}	\includegraphics[width=8cm]{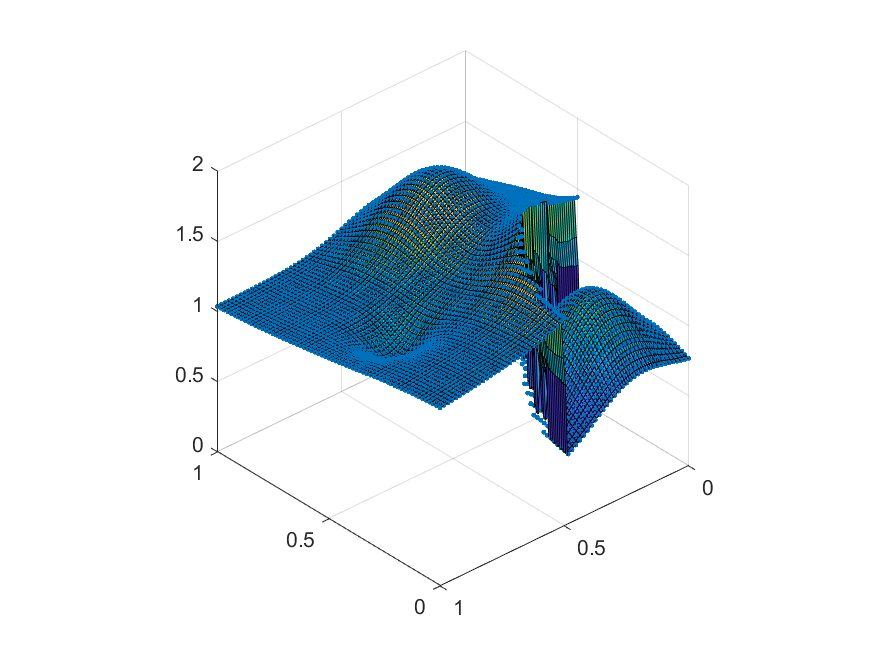}\\
			 \hspace{-1cm}	\includegraphics[width=8cm]{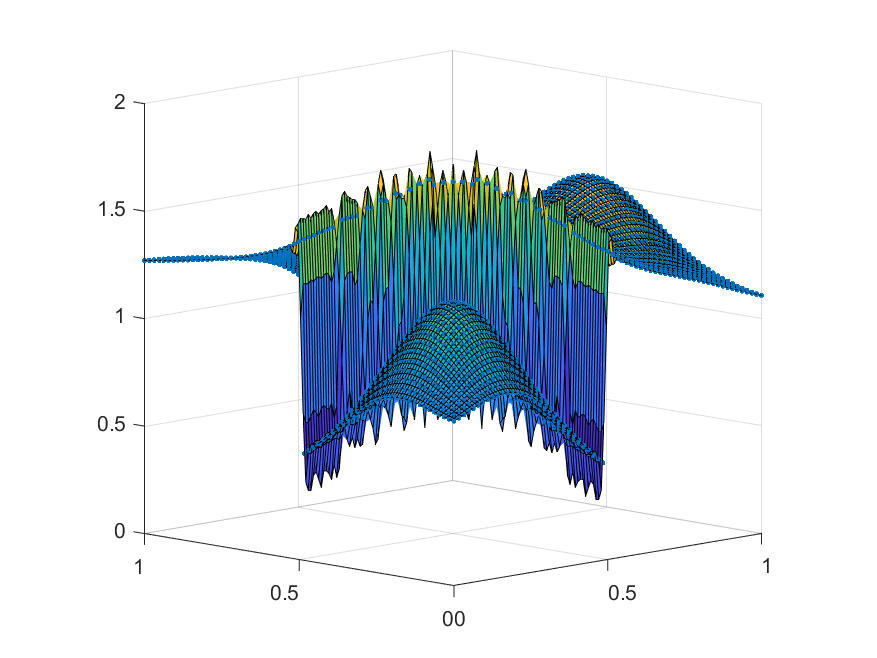} &
		 \hspace{-1cm}	\includegraphics[width=8cm]{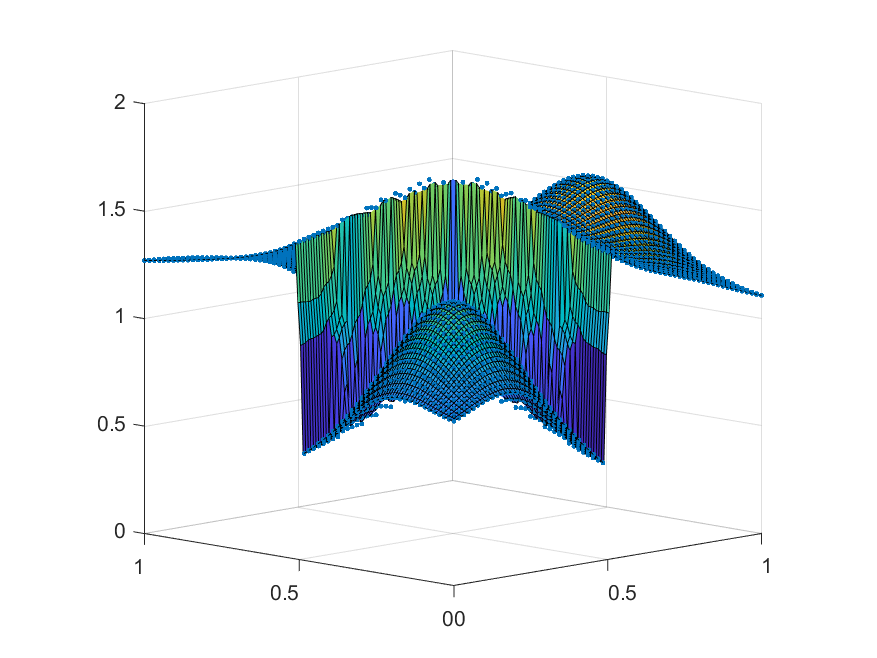}\\
			 \hspace{-1cm}	\includegraphics[width=8cm]{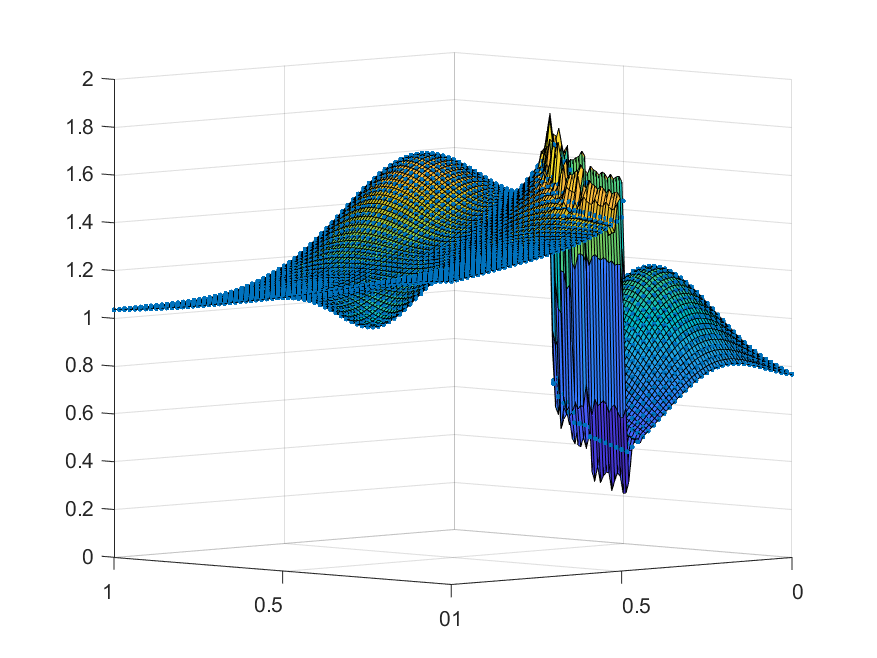} &
		 \hspace{-1cm}	\includegraphics[width=8cm]{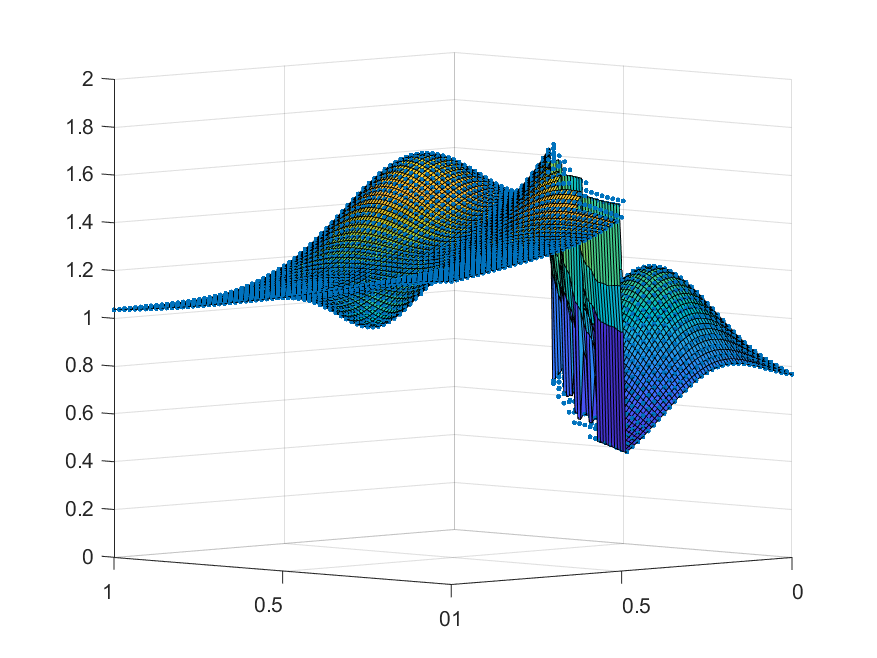}
		\end{tabular}
\end{center}
\caption{Approximation of the function $f_1$ using PUM, and NL-PUM at grid points where we use $\phi_{M_2}$ for the RBF problems and $\phi_{W_2}$ for the partition of unity. Second and third rows are rotations of the plots presented in the first row.}
	\label{figuref1}
	\end{figure}

Let us start with the example presented in the introduction, Sec. \ref{introduction}, where the discontinuity curve is the following
$$\zeta_1(x,y)=x^2+y^2-0.5^2.$$
In this case, we plot in Fig. \ref{figuref1} the result of the experiment with $L=120^2$ equidistant points in the square $[0,1]^2$. It is clear, Fig. \ref{figuref1} right column, that NLPUM presents an adequate behaviour close to the discontinuities whereas the PUM produces some oscillations. When the figure is rotated (Fig. \ref{figuref1} rows 2 and 3), the differences become even more evident.

We also discuss the results when the discontinuities are horizontal, vertical and diagonal lines. For this, we computed a first example with:
 $$\Gamma_2=(\{0.5\}\times [0.5,1] )\cup ([0,0.5]\times\{0.5\}),$$
and a second with
 $$\zeta_3(x,y)=x+y-1.$$

\begin{figure}[hbpt!]
\begin{center}
		\begin{tabular}{cc}
PUM & NL-PUM\\
			 \hspace{-1cm}	\includegraphics[width=8cm]{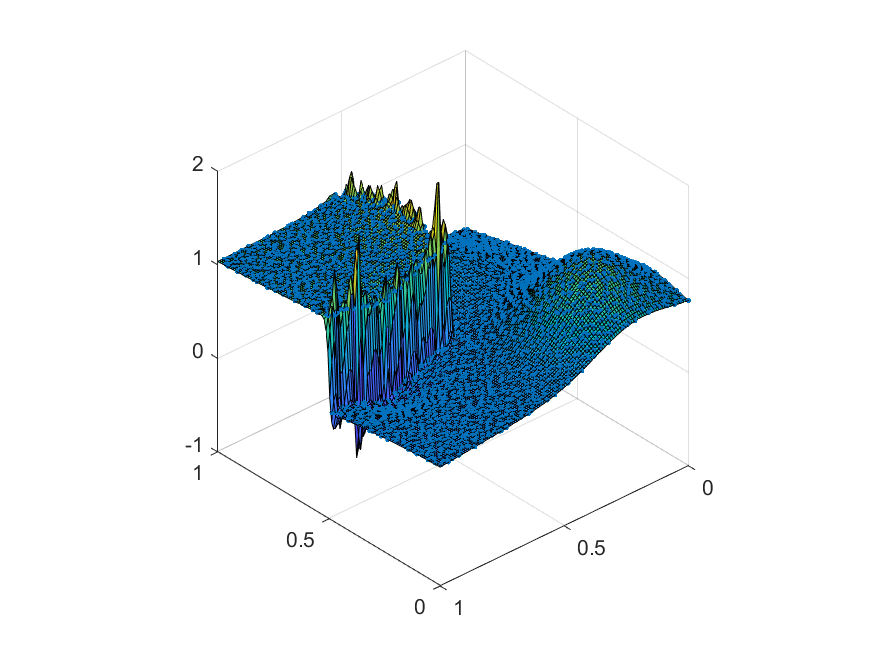} &
		 \hspace{-1cm}	\includegraphics[width=8cm]{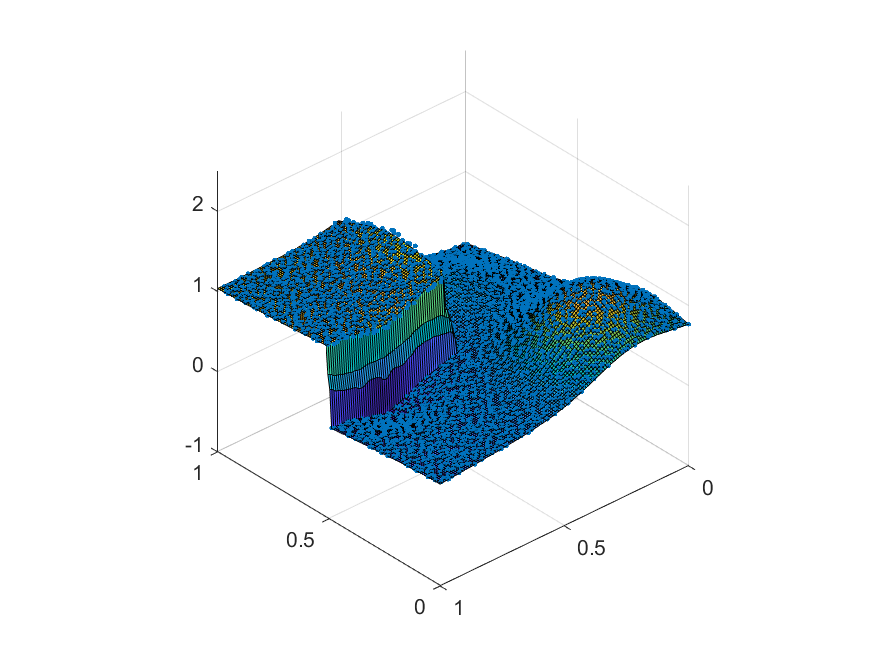}\\
			 \hspace{-1cm}	\includegraphics[width=8cm]{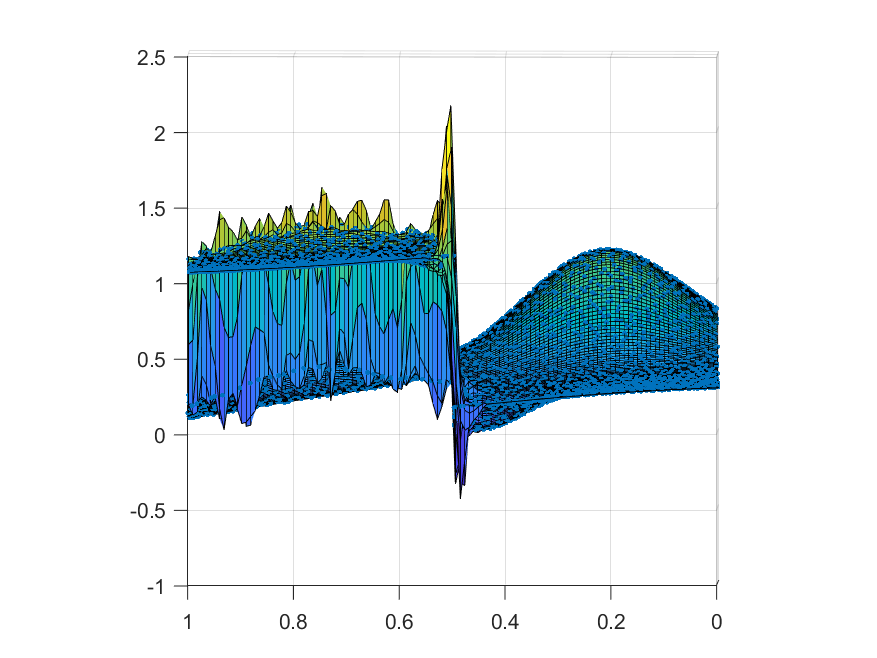} &
		 \hspace{-1cm}	\includegraphics[width=8cm]{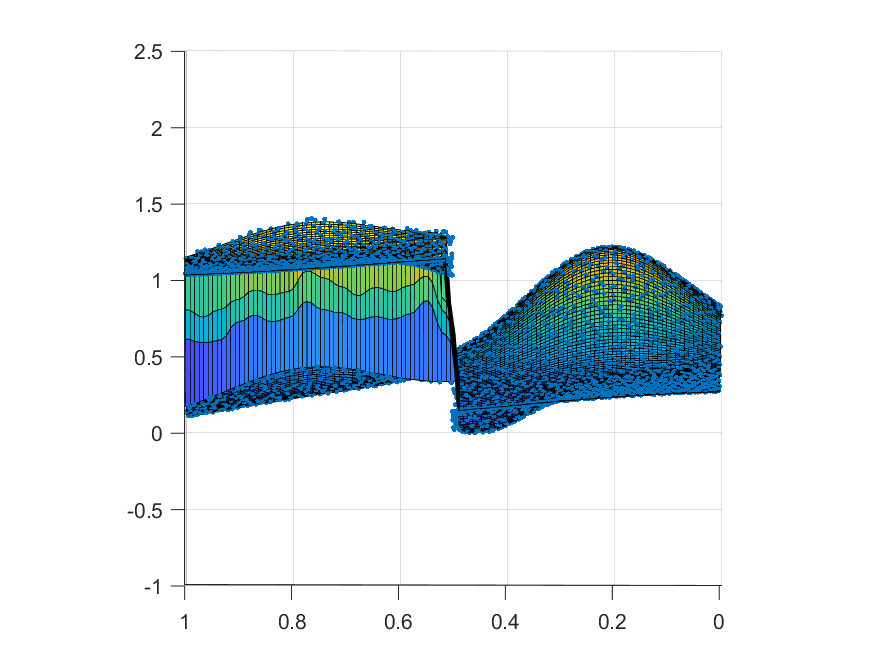}\\
		\end{tabular}
\end{center}
\caption{Approximation to functions $f_2$ using PUM, and NL-PUM at Halton's points using for RBF problems $\phi_{M_2}$ and for partition of unity $\phi_{W_2}$. The second row is a rotation of the plots presented in the first row.}
	\label{figuref2}
	\end{figure}

\begin{figure}[hbpt!]
\begin{center}
		\begin{tabular}{cc}
PUM & NL-PUM\\
 \hspace{-1cm}	\includegraphics[width=8cm]{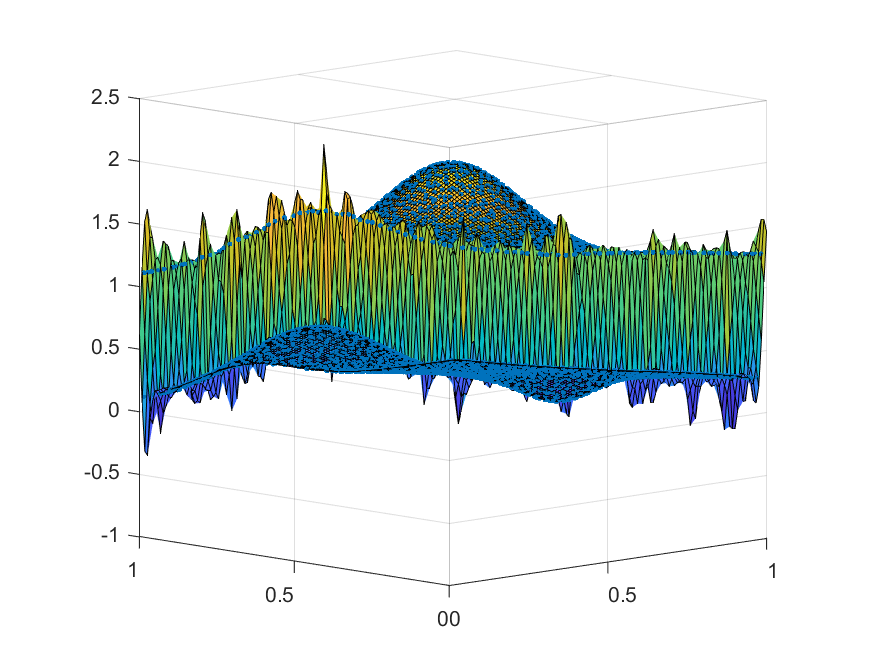} &
		 \hspace{-1cm}	\includegraphics[width=8cm]{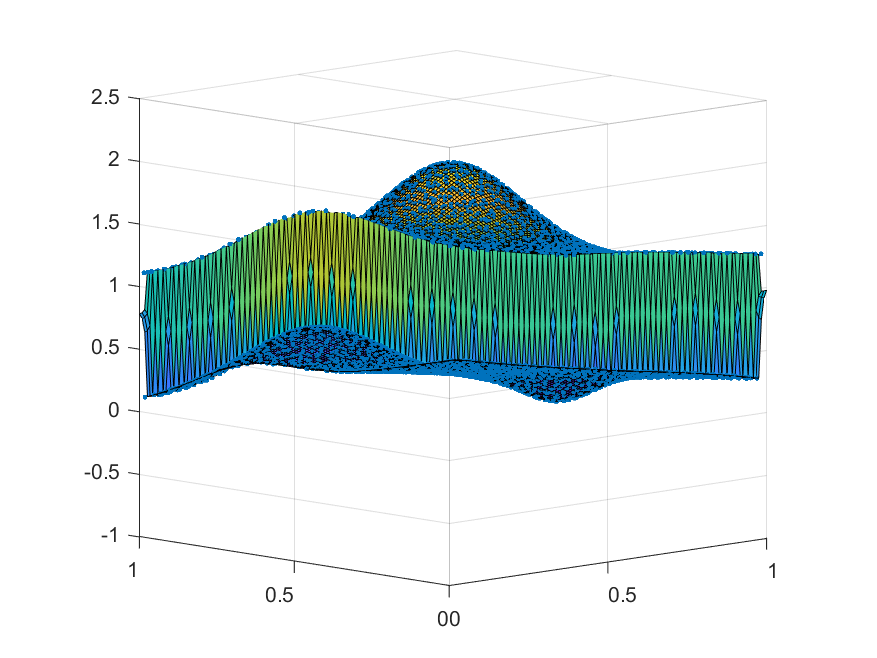}\\
		\end{tabular}
\end{center}
\caption{Approximation to function $f_3$ using PUM, and NL-PUM at Halton's points using for RBF problems $\phi_{M_2}$ and for partition of unity $\phi_{W_2}$.}
	\label{figuref3}
	\end{figure}

 We illustrate the interpolators for the function  $f_2$ in Fig. \ref{figuref2} and for $f_3$ in  Fig. \ref{figuref3}. In these cases we take $N=65^2$ initial Halton's points and we approximate at $L=120^2$ final points. Again, we can see that the new method avoids the oscillations.

\begin{figure}[hbpt!]
\begin{center}
		\begin{tabular}{cc}
PUM & NL-PUM\\
			 \hspace{-1cm}	\includegraphics[width=8cm]{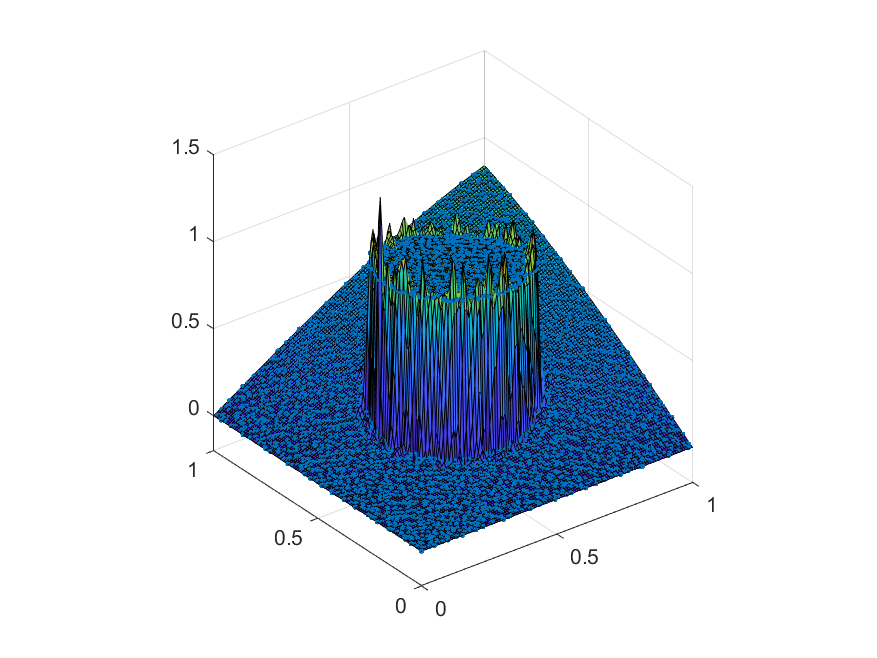} &
		 \hspace{-1cm}	\includegraphics[width=8cm]{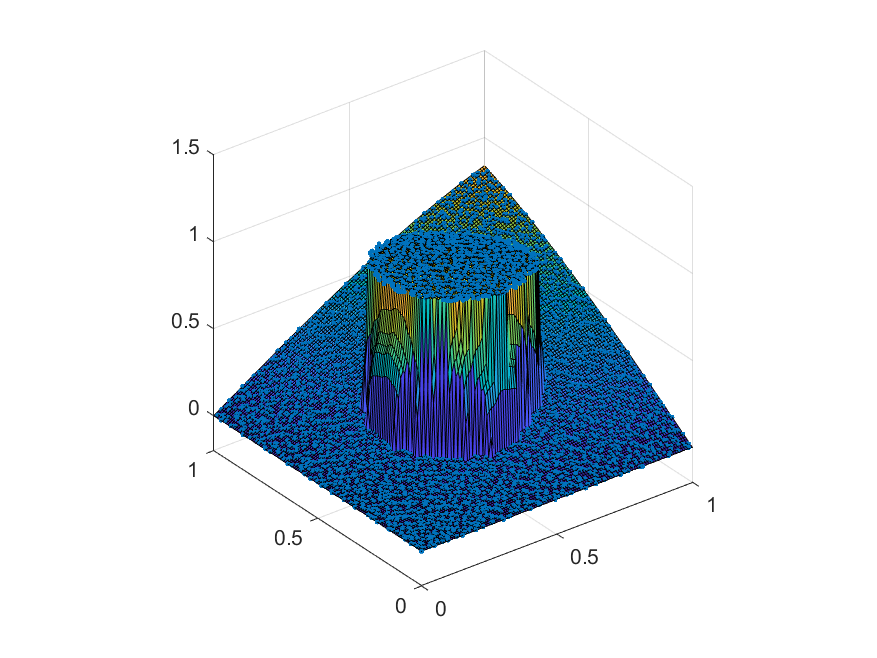}\\
			 \hspace{-1cm}	\includegraphics[width=8cm]{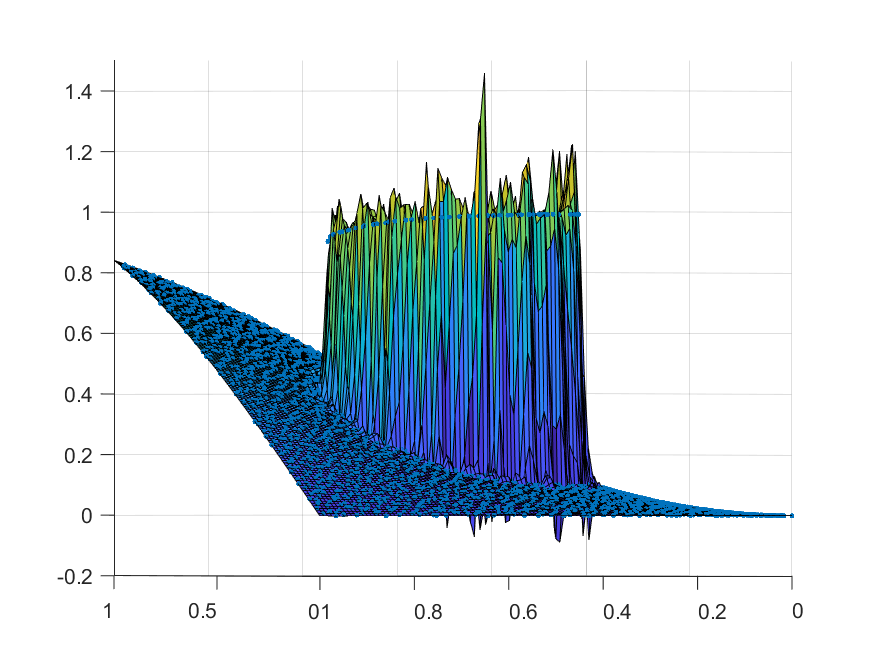} &
		 \hspace{-1cm}	\includegraphics[width=8cm]{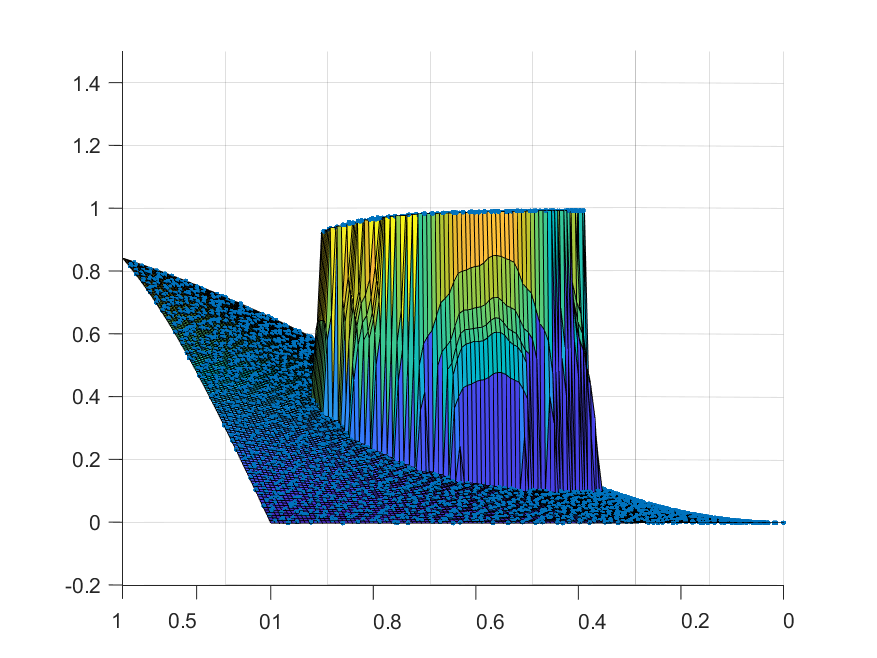}
\end{tabular}
\end{center}
\caption{Approximation to functions $z$ using PUM, and NL-PUM at Halton's points  using for RBF problems $\phi_{M_2}$ and for partition of unity $\phi_{W_2}$.  The second row is a rotation of the plots presented in the first row.}
	\label{figureejemplo2c2}
	\end{figure}

\begin{figure}[hbpt!]
\begin{center}
		\begin{tabular}{cc}
PUM & NL-PUM\\
 \hspace{-1cm}	\includegraphics[width=7cm]{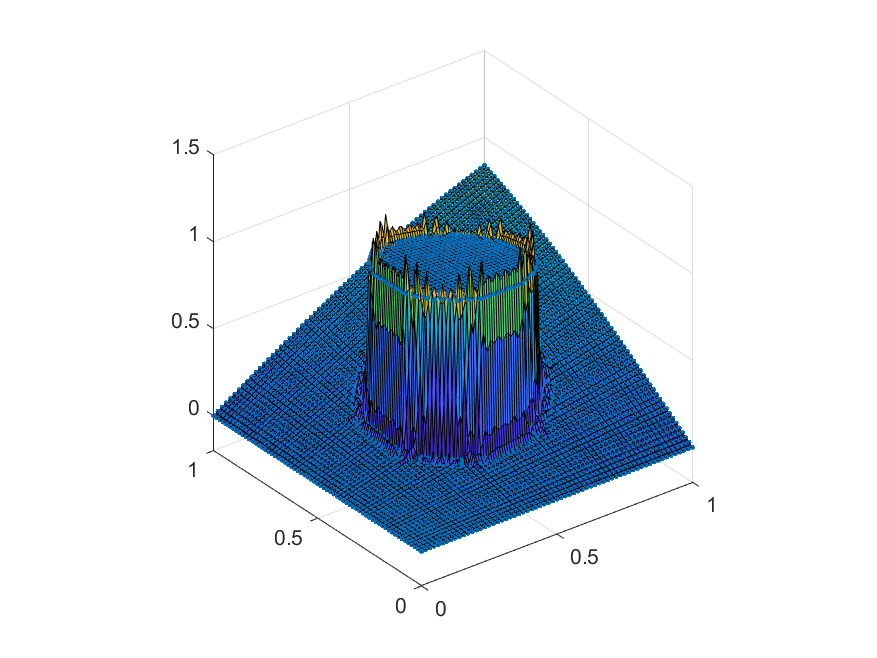} &
		 \hspace{-1cm}	\includegraphics[width=7cm]{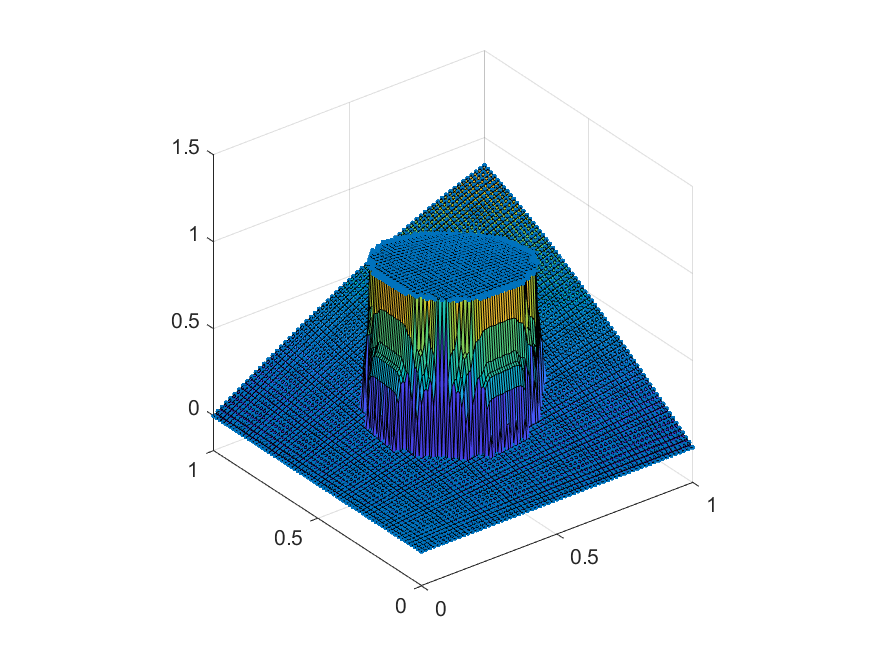}\\
			 \hspace{-1cm}	\includegraphics[width=7cm]{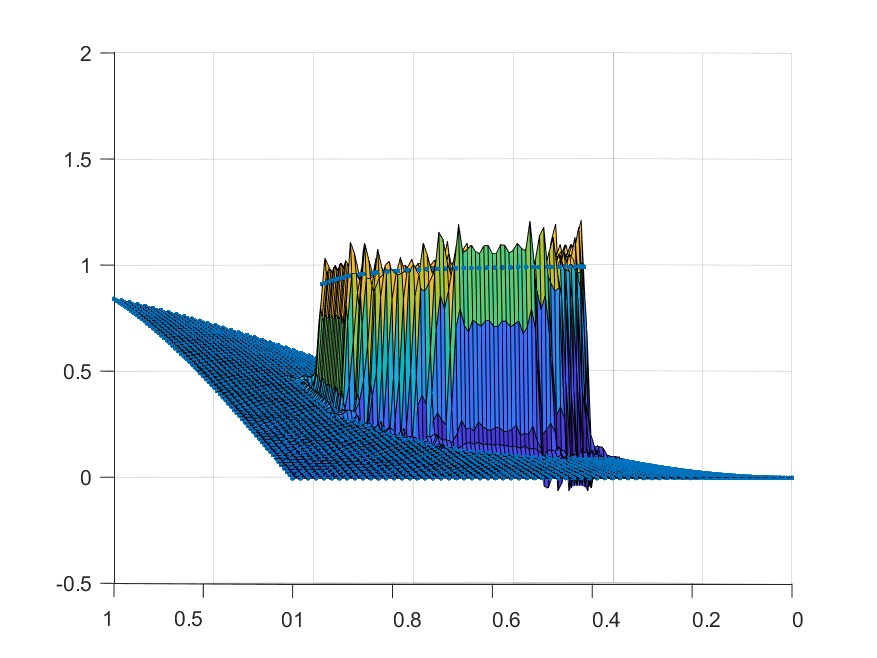} &
		 \hspace{-1cm}	\includegraphics[width=7cm]{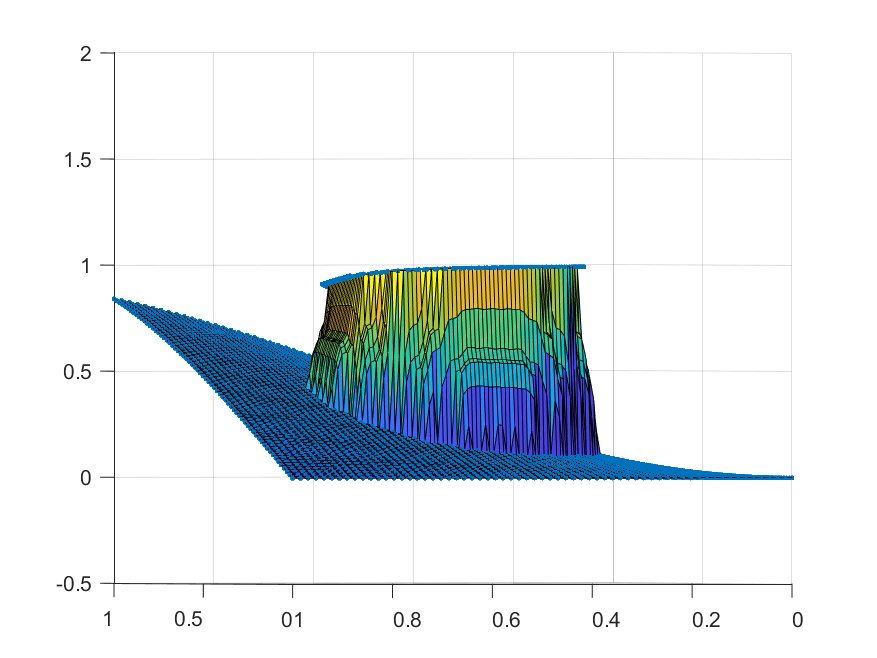}\\
			 \hspace{-1cm}	\includegraphics[width=7cm]{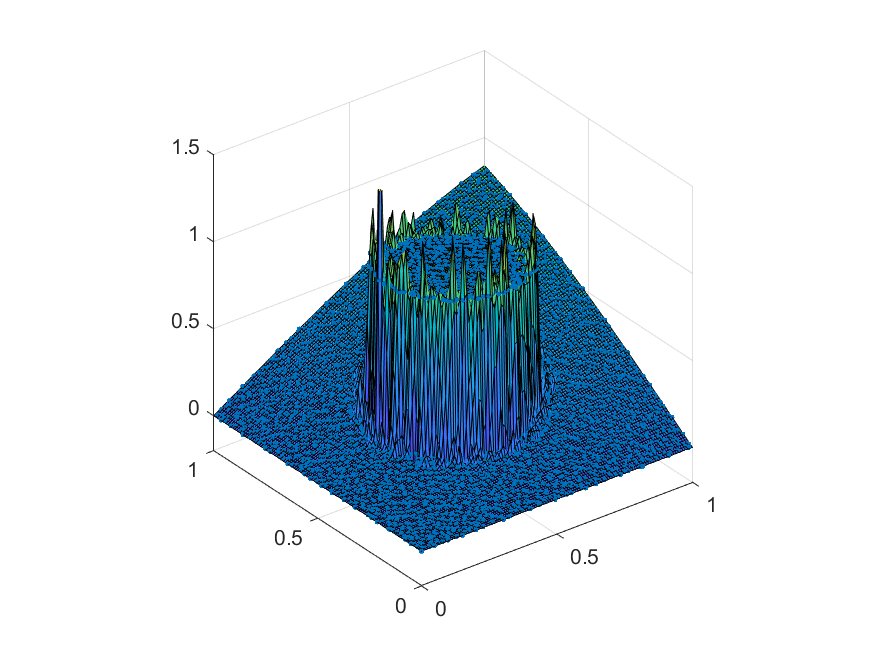} &
		 \hspace{-1cm}	\includegraphics[width=7cm]{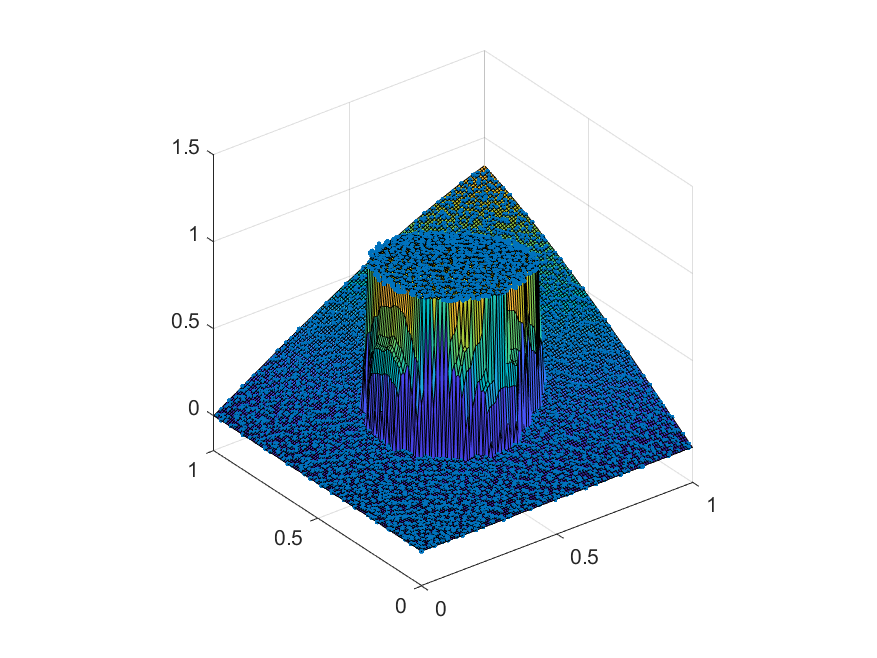}\\
			 \hspace{-1cm}	\includegraphics[width=7cm]{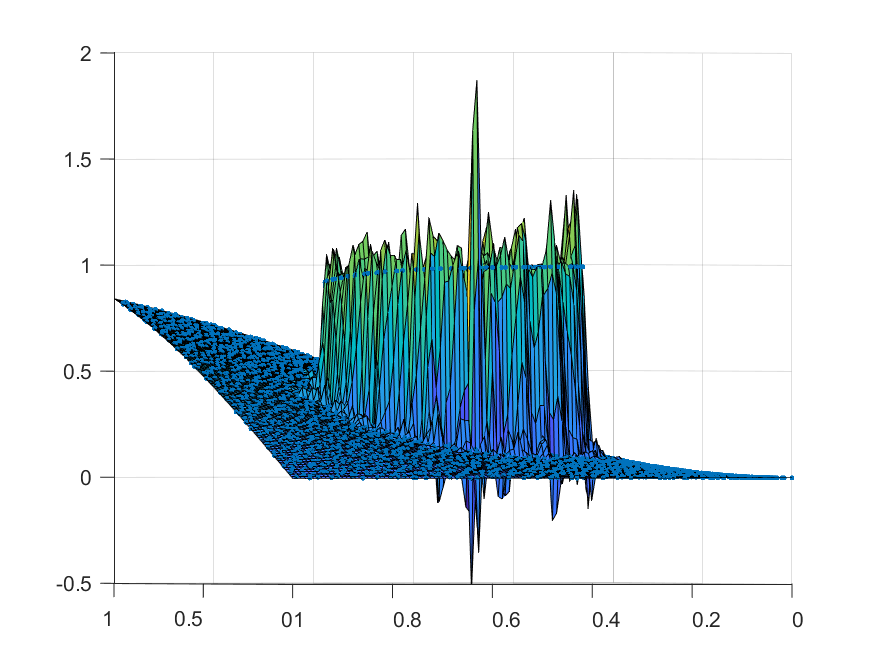} &
		 \hspace{-1cm}	\includegraphics[width=7cm]{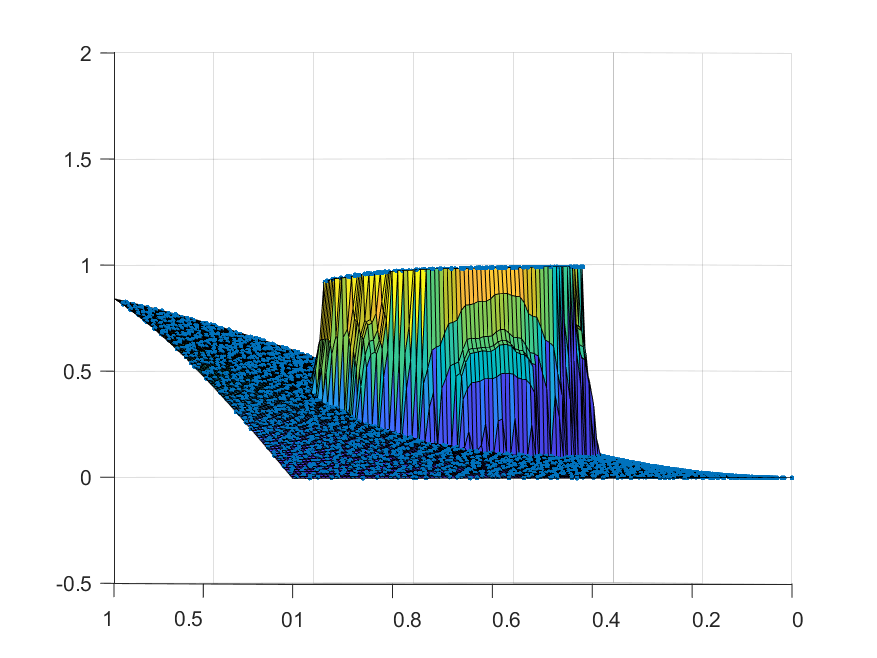}
  		\end{tabular}
\end{center}
\caption{Approximation to function $z$ using PUM, and NL-PUM  where we use  $\phi_{M_4}$  for the RBF problems and $\phi_{W_4}$ for the partition of unity employing gridded data points (two first rows) and Halton's points (two last rows). The second and fourth rows are a rotation of the plots presented in the first row and third row respectively.}
	\label{figureejemplo2c4}
	\end{figure}

To finalize this subsection, we apply the algorithms to the same function employed in \cite{bsplines},
\begin{equation}\label{ejemplonuma}
z(x, y) =
\begin{cases}
\sin(xy), & (x - 0.5)^2 + (y - 0.5)^2 \geq 0.25^2, \\
\cos(xy), & (x - 0.5)^2 + (y - 0.5)^2 < 0.25^2.
\end{cases}
\end{equation}
We plot the result using $\phi_{M_2}$ for RBF problems  and $\phi_{W_2}$ for partition of unity  in Fig. \ref{figureejemplo2c2} using Halton's points. Finally, we illustrate the interpolants employing $\phi_{M_4}$ and $\phi_{W_4}$  in Fig. \ref{figureejemplo2c4} using in the two first rows  gridded data points and the two last rows Halton's point. The interpolators in both cases are very similar. Again, we can see that the new method improves the behaviour close to the discontinuities.
Furthermore, in Fig. \ref{figureejemplo2c4} we perform the experiment for Halton and grided data points.  The oscillations obtained when the classical PUM is applied are highly significant in contrast to the result obtained when the NLPUM is employed.

\section{Conclusions and future work}\label{conc}

The Non-linear Partition of Unity Method (NL-PUM) presented in this study marks a substantial improvement in interpolation techniques, offering a robust solution to challenges posed by discontinuities in datasets. It is founded on integrating Radial Basis Function (RBF) interpolation with Weighted Essentially Non-Oscillatory (WENO) and improves the result obtained by traditional PUM algorithms if the data contains discontinuities. By dynamically adjusting weights through smoothness indicators, NL-PUM selectively conserves the accuracy in smooth regions while minimizing oscillations near discontinuities. These compactly supported, dynamically weighted functions allow for efficient computation without compromising precision.

Numerical experiments confirmed that the NL-PUM maintains the order of accuracy expected from traditional methods in continuous domains while exhibiting superior performance in discontinuous settings. This is achieved by identifying and treating ``contaminated patches'' near discontinuities, effectively eliminating oscillations that commonly arise in such scenarios. The method's ability to consistently produce reliable results, regardless of the orientation or complexity of discontinuities, serves to emphasize its robustness. Additionally, the computational efficiency achieved through advanced data partitioning techniques ensures scalability and applicability across various domains.

While the NL-PUM demonstrates promising capabilities, it also presents certain limitations, such as its reliance on well-tuned parameters. These factors highlight the need for further research to optimize and simplify the method. Future work could focus on automating parameter selection, extending the framework to higher-dimensional problems, and exploring its application in practical contexts such as image processing and fluid dynamics. The potential to integrate NL-PUM with other techniques or enhance its adaptability paves a way for further investigation in the near future. 

%

\end{document}